\documentclass[a4paper,abstracton]{scrartcl}
\usepackage{etex}

\usepackage[utf8]{inputenc}
\usepackage[T1]{fontenc}
\usepackage[english]{babel}
\usepackage{paralist}
\usepackage{amssymb}
\usepackage[ntheorem]{empheq} 
\usepackage[thmmarks,amsmath]{ntheorem}
\usepackage{hyperref}
\usepackage{longtable}
\usepackage{tabu}
\usepackage{jensstuff}
\usepackage{jensntheorem}
\usepackage{lmodern}
\usepackage{framed}
\usepackage{mathrsfs}
\usepackage{dsfont}
\usepackage{listings}
\usepackage{tikz}
\tikzset{>=stealth}
\usetikzlibrary{arrows,positioning}
\usepackage[style=alphabetic,backend=bibtex]{biblatex}
\bibliography{ba}

\definecolor{ao}{rgb}{0.0, 0.5, 0.0}

\allowdisplaybreaks

\title{Exterior powers of F-zips}
\author{Jens Hesse\footnote{University of Paderborn, D-33098 Paderborn, Germany, \texttt{hjens@mail.upb.de}}}
\date{July 28, 2014}

\begin{document}
\maketitle

\begin{abstract}
  An F-zip over a field of positive characteristic is a vector space together with two filtrations whose subquotients are related in a certain way.
  We will define the category of F-zips and some basic constructions in it, especially exterior powers.
  If the ground field is algebraically closed, one can give a classification of F-Zips in terms of combinatorics.
  However,  the way constructions and concepts in the category of F-zips manifest themselves in terms of the classifying invariant, is yet to be fully understood.

  The theory of F-crystals suggests that another invariant might be useful in trying to improve the understanding of F-zips. Given an F-zip, we calculate for every $1$-dimensional F-zip (of which there is essentially one for every integer $d$) $\1(d)$ and every $r\in\Z$ the dimension of the space of F-zip morphisms from $\1(d)$ into the $r$-th exterior power of the given F-zip.

  To make sense of this however, we will have to canonically decompose these spaces of morphisms each into two subspaces that are finite-dimensional over the prime field and its prime field respectively.
  
  One result will then be a way to calculate these numbers for a given isomorphism type.
  Our main result however, is a negative one: The invariant does not classify F-zips.
\end{abstract}

\newpage
\tableofcontents
\newpage

\section{Introduction}
\label{sec:introduction}

When investigating the De Rham cohomology of varieties over a field $k$ of characteristic $p>0$, the notion of an \textit{F-zip} naturally arises. For simplicity, let us assume $k$ to be perfect. Then an F-zip essentially is an $n$-dimensional $k$-vector space together with two filtrations whose subquotients are related in a certain way. In fact, to describe an F-zip up to isomorphism, it is sufficient to say where the dimension jumps in the filtrations and to give a permutation of $1,\dotsc,n$. That permutation will encode the information of how the two filtrations and their respective subquotients are related. The question of when two permutations yield isomorphic F-zips is subtle, but can be  answered (at least if $k$ is algebraically closed). Thus, F-zips can be fully classified. Still, it is not clear, how some information about the F-zip (for example, the property of it being simple or not) can be read off the permutation and how certain operations on F-zips might correspond to certain operations on the permutations associated to them.

This has led to the search for new meaningful invariants of F-zips that ideally even encode the isomorphism class. In this thesis we will examine one such candidate.

More precisely, an F-zip $\mathbf{V}$ is a quadruple consisting of an $n$-dimensional vector space $V$, a descending filtration $C^\bullet$, an ascending  filtration $D_\bullet$ and Frobenius-linear bijections $\varphi_i\colon C^i/C^{i+1}\to D_i/D_{i-1}$ (i.e., $\varphi_i$ is additive and $\varphi_i(ax)=a^p\varphi_i(x)$ for $a\in k$, ${x\in C^i/C^{i+1}}$) and a morphism of F-zips respects the filtrations and makes the obvious squares commute. We will show that there is a natural notion of exterior powers of F-zips (which is just the usual exterior power of the underlying vector space endowed with a suitable F-zip structure as outlined in \cite{1208.3547}). Then one can consider the numbers $\bigl(\dim_{\F_p}\Hom(\1(d),\bigwedge^m\mathbf{V})\bigr)_{(d,m)\in\Z\times\N}$, where the F-zip $\1(d)$ is just $k$ with filtrations $C^\bullet=(\dotsb\supseteq k\supseteq k\supseteq 0\supseteq 0\supseteq \dotsb)$ and $D_\bullet=(\dotsb\subseteq 0\subseteq 0\subseteq k\subseteq k\subseteq \dotsb)$ both jumping in dimension at step $d$, the bijection $C^d/C^{d+1}\to D_d/D_{d-1}$ being simply the Frobenius endomorphism. The hope might be that these provide enough information to maybe even recover the isomorphism class of $\mathbf{V}$.

However, we will point out that these $\Hom$-spaces will in general not even be finite dimensional $\F_p$-vector spaces (also note that is quite clear that they can't be $k$-vector spaces in general). This issue still can be fixed in that we can choose a natural splitting $\Hom(\1(d),\mathbf{V})=U_1\oplus U_2$ where $U_1$ is a finite dimensional $\F_p$-vector space and $U_2$ a finite dimensional $k$-vector space and then consider the dimensions of $U_1$ and $U_2$, but we will also provide an example of a type of F-zips (having height $n=5$) where we can't recover the isomorphism class even after this modification.

In the first section we will define exactly what an F-zip is and what a morphism of F-zips is and we will give some basic examples. In the next section we go on to show that the category of F-zips is $\F_p$-linear (in particular it has finite (co)products) and define the tensor product of F-zips. Things get a bit more complicated when we want to define the exterior power since we will of course want to build it from the tensor product, but taking images, kernels and cokernels of arbitrary morphisms of F-zips will turn out to be problematic; thus making us restrict our attention to morphisms satisfying an extra condition, the \textit{admissible} morphisms. The category of F-zips with said class of admissible morphisms will then be shown to be an exact category (whence the name ``admissible'').

After having discussed all these constructions,  we briefly describe the details of the above-mentioned classification of F-zips using permutations and, as an example, prove  the  full statement in the special case where both filtrations are full (i.e., the dimension always jumps by 0 or 1).

The next section will answer the question of what taking exterior powers means in terms of the associated permutations and the section after  that we will discuss how to still make sense of our original approach of calculating the numbers $\dim_{\F_p}\Hom(\1(d),\mathbf{V})$ by writing $\Hom(\1(d),\mathbf{V})=U_1\oplus U_2$ as described above. We will then be able  to calculate these two dimensions by looking at the permutation associated to $\mathbf{V}$.

Finally, using these results, we will give an example of a  case of two F-zips where all these numbers agree, but their isomorphism classes do not.

\section{F-zips}
\label{sec:f-zips}

In this section we want to introduce the category of F-zips. F-zips are objects of semilinear algebra, so we start with a few basic remarks on that.

\begin{Definition}
  Let $\rho\colon R\to S$ be a ring homomorphism, $M$ an $R$-module and $N$ an $S$-module.

  We then define $\rho^*M:=M\otimes_RS$ (an $S$-module) and note that, via $\rho$, we can consider $N$ also as an $R$-module $\rho_*N$.

  A map $f\colon M\to N$ is called $\rho$-linear if it is $R$-linear considered as a map $M\to \rho_*N$, i.e., if it is additive and $f(rm)=\rho(r)f(m)$ for all $r\in R$ and $m\in M$.
\end{Definition}

\begin{Lemma}
  $\rho^*$ and $\rho_*$ are adjoint functors, i.e.,
  \begin{equation*}
    \Hom_S(\rho^*M,N)\cong\Hom_R(M,\rho_*N)
  \end{equation*}
  functorial in $M$ and $N$.
\end{Lemma}

As our ground field we fix a perfect field $k$ of characteristic $p$. Recall that saying that $k$ is perfect is equivalent to saying that the Frobenius endomorphism
\begin{equation*}
  \sigma\colon k\to k,\; x\mapsto x^p
\end{equation*}
is bijective.

If $V$ is a $k$-vector space, we will write $V^{(p)}$ instead of $\sigma^*V$.

\begin{Definition}
  An \textit{F-zip} over $k$ is a quadruple $\mathbf{V}=(V,C^\bullet,D_\bullet,\varphi_\bullet)$ consisting of a finite dimensional $k$-vector space $V$, a descending chain $C^\bullet=\left(C^i\right)_{i\in\Z}$ of $k$-subspaces of $V$ with $\bigcap_i C^i=0$ ($\Longleftrightarrow\; C^i=0$ for large $i$) and $\bigcup_iC^i=V$ ($\Longleftrightarrow\; C^i=V$ for small $i$), an ascending chain $D_\bullet=\left(D_i\right)_{i\in\Z}$ of $k$-subspaces of $V$ with $\bigcap_iD_i=0$ ($\Longleftrightarrow\; D_i=0$ for  small $i$) and $\bigcup_iD_i=V$ ($\Longleftrightarrow\; D_i=V$ for large $i$) and a family $\varphi_\bullet=\left(\varphi_i\right)_{i\in\Z}$ of $\sigma$-linear bijective maps $\varphi_i\colon C^i/C^{i+1}\to D_i/D_{i-1}$ (alternatively: $\varphi_i\colon (C^i/C^{i+1})^{(p)}\to D_i/D_{i-1}$ $k$-linear isomorphism).

  Notation: $\gr_C^i\mathbf{V}=C^i/C^{i+1}$, $\gr_i^D\mathbf{V}=D_i/D_{i-1}$

  The map $\tau\colon \Z\to\N_0,\; i\mapsto \dim_k\gr_C^i\mathbf{V}=\dim_kC^i-\dim_kC^{i+1}$ is called the \textit{type} of $\mathbf{V}$.
\end{Definition}

\begin{Remark}
  $V$, $\bigoplus_{i\in\Z} \gr_C^i\mathbf{V}$ and $\bigoplus_{i\in\Z} \gr_i^D\mathbf{V}$ all have the same $k$-dimension, i.e., these are isomorphic $k$-vector spaces.
\end{Remark}

\begin{Proof}
  $\dim_k\bigoplus_{i\in\Z} \gr_C^i\mathbf{V}=\sum_{i\in\Z}(\dim_kC^i-\dim_kC^{i+1})=\dim_kV$ and analogously for $V\cong\bigoplus_{i\in\Z} \gr_i^D\mathbf{V}$.
\end{Proof}

\begin{Example}\label{def:tate}
  For $d\in\Z$ define the \textit{Tate F-zip} $\1(d)$ by $V=k$, $C^i=0$ for $i>d$, $C^i=k$ for $i\leq d$, $D_i=0$ for $i<d$ and $D_i=k$ for $i\geq d$ and $\varphi_i=0$ for $i\ne d$ and $\varphi_d=\sigma$.
\end{Example}

The Tate F-zips not only are the most simple example of F-zips (apart from the trivial F-zip), but -- as we shall see later on  -- they also play a special role (as hinted at by the notation $\1(d)$) insofar as the natural notion of tensor products of F-zips is concerned.

\subsection{Morphisms of F-zips}
\label{sec:morphisms-f-zips}

We now have defined F-zips as objects and  complete the definition as follows:

\begin{Definition} \label{def-morph-f-zip}
  A morphism $f\colon \mathbf{V}\to\mathbf{W}$ of F-zips over $k$ is a $k$-linear map $f\colon V\to W$ satisfying the following conditions:
  \begin{enumerate}[(i)]
  \item $f$ respects both the filtrations $C^\bullet$ and $D_\bullet$ of $\mathbf{V}$ and $\mathbf{W}$ in that $f(C^i\mathbf{V})\subseteq C^i\mathbf{W}$ and $f(D_i\mathbf{V})\subseteq D_i\mathbf{W}$ for all $i$.
  \item The following diagram commutes for all $i$:
  \begin{center}
      \begin{tikzpicture}[node distance=4cm, auto]
        \node (grCV) at (0,3) {$\gr_C^i\mathbf{V}$}; 
        \node (grCW) at (4,3) {$\gr_C^i\mathbf{W}$};
        \node (grDW) at (4,0) {$\gr_i^D\mathbf{W}$};
        \node (grDV) at (0,0) {$\gr_i^D\mathbf{V}$};
    
        \draw[->] (grCV) to node {$\gr_C^if$} (grCW);
        \draw[->] (grDV) to node {$\gr_i^Df$} (grDW);
        \draw[->] (grCV) to node[swap] {$\varphi_i(\mathbf{V})$} (grDV);
        \draw[->] (grCW) to node {$\varphi_i(\mathbf{W})$} (grDW);
      \end{tikzpicture}
    \end{center}

    Here $\gr_C^if$ and $\gr_i^Df$ are the natural maps induced by $f$.
  \end{enumerate}
\end{Definition}

We will now take a closer look at the Hom-sets and give some simple examples.

\begin{Example}
  Let $f\colon k\to V$ be a $k$-linear map. Then $f\in \Hom_{(\text{F-zips/}k)}(\1(d),\mathbf{V})$ if and  only if $f(k)\subseteq C^d\mathbf{V}\cap D_d\mathbf{V}$ and the following diagram  commutes:
  \begin{center}
      \begin{tikzpicture}[node distance=4cm, auto]
        \node (grC1) at (0,3) {$k$}; 
        \node (grCV) at (4,3) {$\gr_C^d\mathbf{V}$};
        \node (grDV) at (4,0) {$\gr_d^D\mathbf{V}$};
        \node (grD1) at (0,0) {$k$};
    
        \draw[->] (grC1) to node {$\gr_C^df$} (grCV);
        \draw[->] (grD1) to node {$\gr_d^Df$} (grDV);
        \draw[->] (grC1) to node[swap] {$\sigma$} (grD1);
        \draw[->] (grCV) to node {$\varphi_d(\mathbf{V})$} (grDV);
      \end{tikzpicture}
    \end{center}

    Put differently, $f\in \Hom_{(\text{F-zips/}k)}(\1(d),\mathbf{V})$ if and  only if $f(1)\in C^d\mathbf{V}\cap D_d\mathbf{V}$ and $f(1)+D_{d-1}\mathbf{V}=\varphi_d(\mathbf{V})(f(1)+C^{d+1}\mathbf{V})$.
\end{Example}

\begin{Example}
  Let $f\colon k\to k$ be a $k$-linear map. Then $f\in \Hom_{(\text{F-zips/}k)}(\1(d),\1(d))$ if and only if the following diagram commutes:
  \begin{center}
      \begin{tikzpicture}[node distance=4cm, auto]
        \node (grC1) at (0,2) {$k$}; 
        \node (grCV) at (2,2) {$k$};
        \node (grDV) at (2,0) {$k$};
        \node (grD1) at (0,0) {$k$};
    
        \draw[->] (grC1) to node {$f$} (grCV);
        \draw[->] (grD1) to node {$f$} (grDV);
        \draw[->] (grC1) to node[swap] {$\sigma$} (grD1);
        \draw[->] (grCV) to node {$\sigma$} (grDV);
      \end{tikzpicture}
    \end{center}

    So we have to have $f(1) = \sigma(f(1))$, i.e., $f(1)\in\F_p$.

    Thus $\Hom_{(\text{F-zips/}k)}(\1(d),\1(d))=\{f\in\End_k(k)\suchthat f(1)\in \F_p\}\cong\F_p$.
\end{Example}

\begin{Example}
  Let $d< e$ be integers.

  Then $C^{e}(\1(d))=0$ and $C^e(\1(e))=k$, so that there exists no $f\in\End_k(k)$ with $f\ne 0$ and $f(C^e(\1(e)))\subseteq C^{e}(\1(d))$. Hence, $\Hom_{(\text{F-zips/}k)}(\1(e),\1(d))=0$.

  Similarly, $D_{d}(\1(d))=k$ and $D_d(\1(e))=0$, so that there exists no $f\in\End_k(k)$ with $f\ne 0$ and $f(D_d(\1(d)))\subseteq D_{d}(\1(e))$. Hence, $\Hom_{(\text{F-zips/}k)}(\1(d),\1(e))=0$.
\end{Example}

\begin{Remark}\label{rem:hom-vec}
  Let $a\in k$, $f\in \Hom_{(\text{F-zips/}k)}(\mathbf{V},\mathbf{W})$.

  Then
  \begin{align*}
    \varphi_i(\mathbf{W})\circ (a\,\gr_C^if)&=\sigma(a)\,(\varphi_i(\mathbf{W})\circ\gr_C^i f) \\
    &= \sigma(a)\,(\gr_i^Df\circ\varphi_i(\mathbf{V}))=(\sigma(a)\,\gr_i^Df)\circ\varphi_i(\mathbf{V})
  \end{align*}

  Hence:
  \begin{align*}
    &af\in\Hom_{\text{F-zips/}k}(\mathbf{V},\mathbf{W}) \\
    & \iff a = \sigma(a) \text{ or }\gr_i^Df=\gr_C^if=0 \text{ for all } i \\
    &\iff a\in\F_p \text{ or } \forall i\colon f(D_i\mathbf{V})\subseteq D_{i-1}\mathbf{W} \text{ and } f(C^i\mathbf{V})\subseteq C^{i+1}\mathbf{W}.
  \end{align*}

  It follows that $\Hom_{(\text{F-zips/}k)}(\mathbf{V},\mathbf{W})$ is a $\F_p$-vector space, but not a $k$-vector space  in general.
\end{Remark}

Also, it  will in general not be a finite dimensional $\F_p$-vector space, as illustrated by the following example:

\begin{Example}
  Consider the F-zip $\mathbf{V}$ with $V=k^2=\langle e_1,e_2\rangle$ and $C^{-1}=\langle e_1,e_2\rangle\supseteq C^0=\langle e_1\rangle \supseteq C^1=\langle e_1\rangle \supseteq C^2=0$ and $D_{-2}=0\subseteq D_{-1} = \langle e_1 \rangle \subseteq D_0 = \langle e_1 \rangle \subseteq D_1 = \langle e_1, e_2 \rangle$.

  Then $\tau(i)=0$ for $i\notin\{-1,1\}$ and $\tau(\pm 1)=1$.

  We have $\Hom_{(\text{F-zips/}k)}(\1(0),\mathbf{V})=\{af\suchthat a\in k\}\cong k$ with
  \begin{equation*}
    f\colon k\to \langle e_1,e_2\rangle,\; 1\mapsto e_1.
  \end{equation*}
\end{Example}

For concrete computations involving filtered vector spaces (e.g. F-zips), the following terminology will often be useful:

\begin{Definition} 
Let $K$ be a field and $(V,F^\bullet)$ an $n$-dimensional vector space together with a (not necessarily complete) flag $F^\bullet$ on it.

We shall say that a basis $\{e_i\}_{1\leq i\leq n}$ of $V$ is \textit{adapted} to $F^\bullet$ (or that it is a basis of $(V,F^\bullet)$), if for all $r\in \Z$ there exists $I\subseteq\{1,\dotsc,n\}$ such that $F^r$ is generated by $\{e_i \suchthat i\in I\}.$
\end{Definition}

\begin{Remark}
  There always exists a basis adapted to $F^\bullet,$ as can be seen by successively extending a basis.
\end{Remark}

\section{F-zip constructions}
\label{sec:konstruktionen-mit-f-zips}

We will now investigate the category of F-zips further. We will show that it is an exact $\F_p$-linear tensor category and describe the construction of exterior powers.

Let $\mathbf{U},\mathbf{V},\mathbf{W},\mathbf{Z}$ be F-zips over $k$.

\subsection{Direct sum, (F-zips\texorpdfstring{$/k$}{/k}) is \texorpdfstring{$\F_p$}{F_p}-linear}
\label{sec:direkte-summe}

\begin{Definition} \label{koprodukt-f-zips}
  Define
  \begin{equation*}
    \mathbf{V}\oplus\mathbf{W}=(V\oplus W, \left(C^i\mathbf{V}\oplus C^i\mathbf{W}\right)_i,\left(D_i\mathbf{V}\oplus D_i\mathbf{W}\right)_i,\left(\varphi_i\mathbf{V}\oplus\varphi_i\mathbf{W}\right)_i).
  \end{equation*}

  Then
  \begin{align*}
    \gr_C^i(\mathbf{V}\oplus\mathbf{W})
    &=\left(C^i\mathbf{V}\oplus C^i\mathbf{W}\right)/\left(C^{i+1}\mathbf{V}\oplus C^{i+1}\mathbf{W}\right) \\
    &\cong\left(C^i\mathbf{V}/C^{i+1}\mathbf{V}\right)\oplus\left(C^{i}\mathbf{W}/C^{i+1}\mathbf{W}\right) \\
    &= \gr_C^i\mathbf{V}\oplus\gr_C^i\mathbf{W}
  \end{align*}
  and analogously, $\gr_i^D(\mathbf{V}\oplus\mathbf{W})= \gr_i^D\mathbf{V}\oplus\gr_i^D\mathbf{W}$.
\end{Definition}

\begin{Lemma}\label{lem:addit}
  (F-zips$/k$) is an (additive and) $\F_p$-linear category and the preceding definition indeed gives a description of the (co)product of two F-zips over $k$.
\end{Lemma}

\begin{Proof}
  In Remark \ref{rem:hom-vec} we observed that $\Hom_{(\text{F-zips}/k)}(\mathbf{V},\mathbf{W})\subseteq\Hom_{(k\text{-Mod})}(V,W)$ is an $\F_p$-subspace. Composition then obviously is $\F_p$-bilinear. Thus (F-zips$/k$) is $\F_p$-linear and what remains to be shown is that $\mathbf{V}\oplus\mathbf{W}$ from the preceding definition is indeed the coproduct of $\mathbf{V}$ and $\mathbf{W}$.

To this end, let $\iota_V\colon V\to V\oplus W$ and $\iota_W\colon W\to V\oplus W$ be the natural $k$-linear maps. We have $\iota_V(C^i\mathbf{V})=C^i\mathbf{V}\oplus 0\subseteq C^i\mathbf{V}\oplus C^i\mathbf{W}=C^i(\mathbf{V}\oplus\mathbf{W})$ and analogously ${\iota_V(D_i\mathbf{V})\subseteq D_i(\mathbf{V}\oplus\mathbf{W})}$. Also, ${(\varphi_i\mathbf{V}\oplus\varphi_i\mathbf{W})\circ\gr_C^i\iota_V=\gr_i^D\iota_V\circ\varphi_i(\mathbf{V})\colon\gr_C^i\mathbf{V}\to\gr_i^D\mathbf{V}\oplus\gr_i^D\mathbf{W}}$ since both are equal to $(\varphi_i\mathbf{V}, 0).$ Hence, $\iota_V$ is a morphism of F-zips, and so is $\iota_W$.

Now, given some F-zip $\mathbf{T}$ and two morphisms $f\colon \mathbf{V}\to \mathbf{T}$, $g\colon \mathbf{W}\to \mathbf{T}$, there exists a unique $k$-linear map $u\colon V\oplus W\to T$ with $u\circ\iota_V=f$ and $u\circ\iota_W=g$. We have $u(C^i(\mathbf{V\oplus W}))\subseteq C^i\mathbf{T}$ since $C^i(\mathbf{V\oplus W})=C^i(\mathbf{V})\oplus C^i(\mathbf{W})$ is generated by the elements of $C^i(\mathbf{V})\cong \iota_V(C^i(\mathbf{V}))$ and $C^i(\mathbf{W})\cong \iota_W(C^i(\mathbf{W}))$. Likewise, the diagram
  \begin{center}
      \begin{tikzpicture}[node distance=4cm, auto]
        \node (grCV) at (0,3) {$\gr_C^i\mathbf{V}\oplus\gr_C^i\mathbf{W}$}; 
        \node (grCW) at (7,3) {$\gr_C^i\mathbf{T}$};
        \node (grDW) at (7,0) {$\gr_i^D\mathbf{T}$};
        \node (grDV) at (0,0) {$\gr_i^D\mathbf{V}\oplus \gr_i^D\mathbf{W}$};
    
        \draw[->] (grCV) to node {$\gr_C^iu=\gr_C^if\oplus\gr_C^ig$} (grCW);
        \draw[->] (grDV) to node {$\gr_i^Du=\gr_i^Df\oplus\gr_i^Dg$} (grDW);
        \draw[->] (grCV) to node[swap] {$\varphi_i(\mathbf{V})\oplus\varphi_i(\mathbf{W})$} (grDV);
        \draw[->] (grCW) to node {$\varphi_i(\mathbf{T})$} (grDW);
      \end{tikzpicture}
    \end{center}
    whose commutativity we need to show, can effectively be split into two diagrams that we know to be commutative. Consequently, $u$ is a morphism of F-zips and the lemma is proved.
\end{Proof}

\subsection{Tensor product}
\label{sec:tensorprodukt}

\begin{Definition}
  Define
  \begin{equation*}
    \mathbf{V}\otimes\mathbf{W}=\Bigl(V\otimes W, \Bigl(\sum_{i+j=r}C^i\mathbf{V}\otimes C^j\mathbf{W}\Bigr)_r, \Bigl(\sum_{i+j=r}D_i\mathbf{V}\otimes D_j\mathbf{W}\Bigr)_r, \Bigl(\bigoplus_{i+j=r}\varphi_i\mathbf{V}\otimes\varphi_j\mathbf{W}\Bigr)_r\Bigl).
  \end{equation*}
\end{Definition}

\begin{Lemma}
  (F-zips$/k$) is a tensor category, the unit object being $\1(0)$ (defined in Example \ref{def:tate}), in the sense that there are natural isomorphisms ${(\mathbf{V}\otimes\mathbf{W})\otimes\mathbf{Z}\cong\mathbf{V}\otimes(\mathbf{W}\otimes\mathbf{Z})}$, ${\1(0)\otimes\mathbf{V}\cong\mathbf{V}}$ and ${\mathbf{V}\otimes\1(0)\cong\mathbf{V}}$, which are compatible in the obvious ways.
\end{Lemma}

\begin{Proof}
  ($k$-Mod) is a tensor category and one just has to check that the isomorphisms $(V\otimes W)\otimes Z\cong V\otimes(W\otimes Z)$, \ldots are compatible with the filtrations and make the  relevant squares commute, which is immediate.
\end{Proof}

\begin{Lemma}
  The natural isomorphism $(U\oplus V)\otimes W\to (U\otimes W)\oplus (V\otimes W)$ is an isomorphism of F-zips.
\end{Lemma}

\begin{Proof}
  We observe  that
  \begin{align*}
    C^r((\mathbf{U}\oplus\mathbf{V})\otimes\mathbf{W})
    &= \sum_{i+j=r}(C^i\mathbf{U}\oplus C^i\mathbf{V})\otimes C^j\mathbf{W} \\
    &\cong \sum_{i+j=r}\bigl((C^i\mathbf{U}\otimes C^j\mathbf{W}) \oplus (C^i\mathbf{V}\otimes C^j\mathbf{W})\bigr) \\
    &= \sum_{i+j=r}(C^i\mathbf{U}\otimes C^j\mathbf{W}) \oplus  \sum_{i+j=r}(C^i\mathbf{V}\otimes C^j\mathbf{W}) \\
    &= C^r(\mathbf{U}\otimes \mathbf{W})\oplus C^r(\mathbf{V}\otimes \mathbf{W}).
  \end{align*}

  Like  this it follows that the natural isomorphism respects the filtration, i.e., (i) from Definition \ref{def-morph-f-zip} is satisfied. It is then straightforward to verify that (ii) is satisfied as well.
\end{Proof}

\begin{Lemma}
  There are natural isomorphisms
  \begin{align*}
    \gr_C^r(\mathbf{V}\otimes\mathbf{W})\cong\bigoplus_{i+j=r}\gr_C^i\mathbf{V}\otimes\gr_C^j\mathbf{W},\\
    \gr_r^D(\mathbf{V}\otimes\mathbf{W})\cong\bigoplus_{i+j=r}\gr_i^D\mathbf{V}\otimes\gr_j^D\mathbf{W}.
  \end{align*}
\end{Lemma}

\begin{Proof}  (of the first equation, the second is proved in the same way)

  Here we consider a finite dimensional vector space $V$ endowed with a descending filtration $C^\bullet$.

  The finite dimensional vector spaces  endowed with a descending filtration form a category (Fil$^\bullet$-Vec) in the obvious way (cf. Definition \ref{def-morph-f-zip}) that has finite coproducts (cf. Definition \ref{koprodukt-f-zips}).

  One may write $(V,C^\bullet)=(V_1,C_1^\bullet)\oplus\dotsb\oplus(V_n,C_n^\bullet)$ with $\forall i\colon\dim V_i=1$: Choose a basis $e_1,\dotsc,e_n$ of $V$ with $C^i=\langle e_j \suchthat h(i)\leq j\leq n\rangle$, where $h\colon \Z\to \{1,\dotsc,n\}\cup\{\infty\}$ is monotonically increasing. Then set $V_\nu=\langle e_\nu\rangle$ and $C_\nu^i=V_\nu$ if $i\leq s_\nu:=\sup\{i\suchthat h(i)\leq \nu\}$ and $C_\nu^i=0$ else.

  Now
  \begin{equation*}
    C^r(\mathbf{V}_\nu\otimes \mathbf{W})
    = \sum_{j\geq r- {s_\nu}}V_\nu\otimes C^j(\mathbf{W})
    = V_\nu\otimes C^{r-{s_\nu}}(\mathbf{W}),
  \end{equation*}
  hence
  \begin{align*}
    \gr_C^r(\mathbf{V}_\nu\otimes \mathbf{W}) 
    &= C^r(\mathbf{V}_\nu\otimes \mathbf{W})/C^{r+1}(\mathbf{V}_\nu\otimes \mathbf{W}) \\
    &\cong (V_\nu\otimes C^{r-{s_\nu}}(\mathbf{W}))/(V_\nu\otimes C^{r-{s_\nu}+1}(\mathbf{W})) \\
    &\cong V_\nu\otimes(C^{r-{s_\nu}}(\mathbf{W})/C^{r-{s_\nu}+1}(\mathbf{W})) \\
    &\cong \gr_C^{s_\nu}(\mathbf{V}_\nu)\otimes \gr_C^{r-{s_\nu}}(\mathbf{W}) \\
    &\cong \bigoplus_{i+j=r}\Bigl(\gr_C^i(\mathbf{V}_\nu)\otimes \gr_C^{j}(\mathbf{W})\Bigr).
  \end{align*}

  Consequently,
  \begin{align*}
    \gr_C^r(\mathbf{V}\otimes\mathbf{W})
    &\cong\gr_C^r(\Bigl(\bigoplus_{\nu=1}^n\mathbf{V}_\nu\Bigr)\otimes\mathbf{W}) \\
    &\cong \bigoplus_{\nu=1}^n\gr_C^r(\mathbf{V}_\nu\otimes\mathbf{W}) \\
    &\cong \bigoplus_{\nu=1}^n\bigoplus_{i+j=r}\bigl(\gr_C^i\mathbf{V}_\nu\otimes\gr_C^j\mathbf{W}\bigr) \\
    &\cong \bigoplus_{i+j=r}\bigoplus_{\nu=1}^n\bigl(\gr_C^i\mathbf{V}_\nu\otimes\gr_C^j\mathbf{W}\bigr) \\
    &\cong \bigoplus_{i+j=r}\Bigl(\bigoplus_{\nu=1}^n\gr_C^i\mathbf{V}_\nu\Bigr)\otimes\gr_C^j\mathbf{W} \\
    &\cong \bigoplus_{i+j=r}\gr_C^i\Bigl(\bigoplus_{\nu=1}^n\mathbf{V}_\nu\Bigr)\otimes\gr_C^j\mathbf{W} \\
    &\cong \bigoplus_{i+j=r}\gr_C^i\mathbf{V}\otimes\gr_C^j\mathbf{W}.
  \end{align*}
\end{Proof}

\begin{Remark}
  In the last equation we are making the following identifications:
  \begin{align*}
    \sum_{i+j=r} \Bigl(\sum_\nu A^\nu_{ij}e_\nu\Bigr)\otimes w_{ij} + C^{r+1}(\mathbf{V}\otimes\mathbf{W})
    &\mathop{\widehat=} \sum_{i+j=r} (A^\nu_{ij}e_\nu)_\nu\otimes w_{ij} + C^{r+1}(\bigoplus_{\nu}\mathbf{V_\nu}\otimes\mathbf{W}) \\
    &\mathop{\widehat=} \left(\sum_{i+j=r} A^\nu_{ij}e_\nu\otimes w_{ij} + C^{r+1}(\mathbf{V}_\nu\otimes\mathbf{W})\right)_\nu \\
    &\mathop{\widehat=} \left(\left((A^\nu_{ij}e_\nu + C^{i+1}\mathbf{V}_\nu)\otimes (w_{ij}+C^{j+1}\mathbf{W})\right)_{i+j=r}\right)_\nu \\
    &\mathop{\widehat=} \left(\left((A^\nu_{ij}e_\nu + C^{i+1}\mathbf{V}_\nu)\otimes (w_{ij}+C^{j+1}\mathbf{W})\right)_\nu\right)_{i+j=r} \\
    &\mathop{\widehat=} \left(\left(A^\nu_{ij}e_\nu + C^{i+1}\mathbf{V}_\nu\right)_\nu\otimes (w_{ij}+C^{j+1}\mathbf{W})\right)_{i+j=r} \\
    &\mathop{\widehat=} \left(\Bigl(\sum_\nu A^\nu_{ij}e_\nu + C^{i+1}\mathbf{V}\Bigr)\otimes (w_{ij}+C^{j+1}\mathbf{W})\right)_{i+j=r},
  \end{align*}
  where $w_{ij}\in C^j\mathbf{W}, A_{ij}^\nu\in k,V_\nu=\langle e_\nu\rangle$.
\end{Remark}

\begin{Remark}\label{rem:gr-tensor}
  Generalizing what was previously shown, one has the following identities:
  \begin{align*}
    C^r(\bigotimes_{\nu=1}^m\mathbf{V}_\nu) &= \sum_{i_1+\dotsb+i_m=r}\bigotimes_{\nu=1}^mC^{i_\nu}\mathbf{V}_\nu, \\
    D_r(\bigotimes_{\nu=1}^m\mathbf{V}_\nu) &= \sum_{i_1+\dotsb+i_m=r}\bigotimes_{\nu=1}^mD_{i_\nu}\mathbf{V}_\nu, \\
    \gr_C^r(\bigotimes_{\nu=1}^m\mathbf{V}_\nu) &= \bigoplus_{i_1+\dotsb+i_m=r}\bigotimes_{\nu=1}^m\gr_C^{i_\nu}\mathbf{V}_\nu, \\
    \gr_r^D(\bigotimes_{\nu=1}^m\mathbf{V}_\nu) &= \bigoplus_{i_1+\dotsb+i_m=r}\bigotimes_{\nu=1}^m\gr_{i_\nu}^D\mathbf{V}_\nu.
  \end{align*}
\end{Remark}

\begin{Remark}
  $\mathbf{V}\otimes\1(d)$ is ``$\mathbf{V}$ shifted to the right by $d$'' in the following sense:
  \begin{equation*}
    C^r(\mathbf{V}\otimes\1(d))=\sum_{i+j=r,\;j\leq d}C^i=C^{r-d}, \quad \text{etc.}
  \end{equation*}
\end{Remark}

\subsection{Admissible morphisms, images, kernels, cokernels}
\label{sec:zulass-morph-bild}

Our next aim is to construct the kernel (and the cokernel and the image) of a morphism of F-zips. To this end, we will want to take its kernel as a $k$-linear map and then endow it with an appropriate F-zip structure. However, this construction will only work for a certain class of nice morphisms. When we restrict our attention to this class of morphisms, we will also be able to make sense of exact sequences in the category of F-zips.

\begin{Definition}\label{def:adm}
  $f\in\Hom_{(\text{F-zips/}k)}(\mathbf{V},\mathbf{W})$ is called \textit{admissible} if $f(C^i\mathbf{V})=f(V)\cap C^i\mathbf{W}$ and $f(D_i\mathbf{V})=f(V)\cap D_i\mathbf{W}$ for all $i$.
\end{Definition}

\begin{Lemma} \label{idempotent}
  Let $f\in\Hom_{(\text{F-zips/}k)}(\mathbf{V},\mathbf{V})$ with $f^2=\lambda f$ for some $\lambda\in k^\times$.

  Then it follows that $f$ is admissible.
\end{Lemma}

\begin{Proof}
  We have to show that $f(C^i)\supseteq f(V)\cap C^i$ and $f(D_i)\supseteq f(V)\cap D_i$.

  Let $x\in f(V)\cap C^i$. Then $x=f(v)$ for some $v\in V$. Thus $f(x)=f^2(v) = \lambda f(v) = \lambda x$, hence $x = f(x/\lambda)\in f(C^i)$.

  The second inclusion is proved analogously.
\end{Proof}

As mentioned before, we can now construct images, kernels and cokernels of admissible morphisms in a natural way. 

\begin{Definition}\label{def:bild}
  Let $f\in\Hom_{(\text{F-zips/}k)}(\mathbf{V},\mathbf{W})$ be admissible. Define
  \begin{align*}
    \im(f):=(f(V),  f(C^\bullet\mathbf{V})=f(V)\cap C^\bullet\mathbf{W}, f(D_\bullet\mathbf{V})=f(V)\cap D_\bullet\mathbf{W}, \left.\varphi_\bullet(\mathbf{W})\right|_{\gr_C^\bullet(\im(f))}).
  \end{align*}

  $\varphi_i(\mathbf{W})$ maps $\gr_C^i(\im(f))$ into $\gr_i^D(\im(f))$, because
  \begin{align*}
    \gr_C^i(\im(f))
    &=(f(V)\cap C^i\mathbf{W})/(f(V)\cap C^{i+1}\mathbf{W}) \\
    &\cong\{\overline{f(c)}=\gr_C^if(\overline{c})\in C^i\mathbf{W}/C^{i+1}\mathbf{W} \suchthat c\in C^i\mathbf{V}\}\subseteq C^i\mathbf{W}/C^{i+1}\mathbf{W}
  \end{align*}
  and for $c\in C^i\mathbf{V}$:
  \begin{align*}
    (\varphi_i(\mathbf{W})\circ \gr_C^if)(\overline{c}) = (\gr_i^Df\circ\varphi_i(\mathbf{V}))(\overline{c})
    &\in \{\overline{f(d)}=\gr_i^Df(\overline{d})\in D_i\mathbf{W}/D_{i-1}\mathbf{W} \suchthat d\in D_i\mathbf{V}\} \\
    &\cong (f(V)\cap D_i\mathbf{W})/(f(V)\cap D_{i-1}\mathbf{W})=\gr_i^D(\im(f)).
  \end{align*}
\end{Definition}

\begin{Definition}
  Let $f\in\Hom_{(\text{F-zips/}k)}(\mathbf{V},\mathbf{W})$ be admissible. Define
  \begin{align*}
    \ker(f):=(\ker(f),  \ker(f)\cap C^\bullet\mathbf{V}, \ker(f)\cap D_\bullet\mathbf{V}, \left.\varphi_\bullet(\mathbf{V})\right|_{\gr_C^\bullet(\ker(f))}).
  \end{align*}

  $\varphi_i(\mathbf{V})$ maps $\gr_C^i(\ker(f))$ into $\gr_i^D(\ker(f))$ because of the following:

  For $d\in D_i\mathbf{V}$, we have
  \begin{align*}
    \gr_i^Df(\overline{d})=0 &\iff f(d)\in f(V)\cap D_{i-1}\mathbf{W}=f(D_{i-1}\mathbf{V}) \\
    &\iff f(d) = f(d') \text{ for some } d'\in D_{i-1}\mathbf{V}
  \end{align*}
  and
  \begin{align*}
    \gr_C^i(\ker(f)) &= (\ker(f)\cap C^i\mathbf{V})/(\ker(f)\cap C^{i+1}\mathbf{V}) \\
    &\cong \{\overline{c}\in C^i/C^{i+1}\suchthat c\in C^i\mathbf{V}, f(c)=0\}\subseteq\gr_C^i\mathbf{V}
  \end{align*}
  and for $c\in C^i\mathbf{V}$ with $f(c)=0$ (i.e., $\overline{c}\in\gr_C^i(\ker(f))$) and $d\in D_i\mathbf{V}$ with $\overline{d}=\varphi_i(\mathbf{V})(\overline{c})$ we have
  \begin{equation*}
    0=(\varphi_i(\mathbf{W})\circ \gr_C^if)(\overline{c}) = (\gr_i^Df\circ\varphi_i(\mathbf{V}))(\overline{c})=\gr_i^Df(\overline{d}),
  \end{equation*}
  so that
  \begin{equation*}
    f(d) = f(d') \text{ for some } d'\in D_{i-1}\mathbf{V}.
  \end{equation*}

  Now
  \begin{equation*}
    f(d-d')=0 \text{ and } \overline{d}=\overline{d-d'} \text{ in } D_i/D_{i-1},
  \end{equation*}
  so that
  \begin{equation*}
    \varphi_i(\mathbf{V})(\overline{c})
    \in\{\overline{d}\in D_i/D_{i+1}\suchthat d\in D_i\mathbf{V}, f(d)=0\}\cong \gr_i^D(\ker(f)).
  \end{equation*}
\end{Definition}

\begin{Definition}
  Let $f\in\Hom_{(\text{F-zips/}k)}(\mathbf{V},\mathbf{W})$ be admissible, $\pi\colon W\to\coker(f)$ the natural projection. Define
  \begin{equation*}
    \coker(f)=(\coker(f), \pi(C^\bullet\mathbf{W}), \pi(D_\bullet\mathbf{W}), \overline{\varphi_\bullet(\mathbf{W})}).
  \end{equation*}

  This makes sense because of the following:

  \begin{equation*}
    \pi(C^i\mathbf{W}) \cong C^i\mathbf{W}/(f(V)\cap C^{i}\mathbf{W})=C^i\mathbf{W}/f(C^i\mathbf{V})
  \end{equation*}
  and
  \begin{align*}
    \gr_C^i(\coker(f))
    &= \pi(C^i\mathbf{W})/\pi(C^{i+1}\mathbf{W})
      = (C^i\mathbf{W}/f(C^i\mathbf{V}))/(C^{i+1}\mathbf{W}/f(C^{i+1}\mathbf{V})) \\
      &\cong C^i\mathbf{W}/(C^{i+1}\mathbf{W} + f(C^{i}\mathbf{V})).
  \end{align*}
\end{Definition}

This completes the discussion of the technicalities. The main point is summarized by the following proposition.

\begin{Proposition}
  Every admissible morphism $f\colon\mathbf{V}\to\mathbf{W}$ of F-zips has a kernel, a cokernel and an image and those are preserved by the forgetful functor into the  category of $k$-vector spaces.

  More precisely: The constructions of $\ker(f)$, $\coker(f)$, $\im(f)$ given above indeed yield a kernel, a cokernel and an image of $f$.

  Moreover, the associated morphisms $\ker(f)\to\mathbf{V}$, $\mathbf{W}\to\coker(f)$ and $\im(f)\to\mathbf{W}$ all are admissible.
\end{Proposition}

The proof of course is straightforward. The only difficult part is seeing that the kernel resp. the cokernel resp. the image of $f$ in the category of F-zips can be endowed with a natural F-zip structure. But  we have already shown that this is the case if $f$ is admissible.
 
Next, we make precise the statement that the  category of F-zips is exact.

The idea is as follows: The building blocks of exact sequences (in, say, an abelian category) are short exact sequences, i.e., pairs $A'\xrightarrow{i}A\xrightarrow{p}A''$ of morphisms where $i$ is a kernel of $p$ and $p$ is a cokernel of $i$. And in fact, all we need to make sense of exact sequences is a class of such kernel-cokernel pairs, which is well behaved in certain ways. In our particular example, we will want the morphisms $i$ and $p$ to be  admissible. In fact, this is where the term ``admissible'' stems from in the first place.

We now cite the definition of an exact category from \cite{0811.1480} (Definition 2.1):

\begin{Definition}
  Let $\mathcal{A}$ be an additive category.

  If a class $\mathscr{E}$ of kernel-cokernel pairs is fixed, an \textit{admissible monic} is a morphism $i$ for which there exists a morphism $p$ such that $(i,p)\in\mathscr{E}$. \textit{Admissible epics} are defined dually.

An \textit{exact structure} on $\mathcal{A}$ is a class $\mathscr{E}$ of kernel–cokernel pairs which is closed under
isomorphisms and satisfies the following axioms:
  \begin{longtabu}{lX}
{[E0]} &  For all objects $A\in\mathcal{A}$,  $\id_A$ is an admissible monic. \\
{[E0$^{\text{op}}$]} & For all objects $A\in\mathcal{A}$,  $\id_A$ is an admissible epic. \\
{[E1]} & The class of admissible monics is closed under composition. \\
{[E1$^{\text{op}}$]} & The class of admissible epics is closed under composition. \\
{[E2]} & The push-out of an admissible monic along an arbitrary morphism exists and yields an admissible monic. \\
{[E2$^{\text{op}}$]} & The pull-back of an admissible epic along an arbitrary morphism exists and yields
an admissible epic.
  \end{longtabu}

An exact category is a pair ($\mathcal{A}$,$\mathscr{E}$) consisting of an additive category $\mathcal{A}$ and an exact
structure $\mathscr{E}$ on $\mathcal{A}$.

Elements of $\mathscr{E}$ are called \textit{short exact sequences}.
\end{Definition}

\begin{Lemma}\label{lem:adm}
  A morphism $\mathbf{V}\to\mathbf{W}$ in (F-zips$/k$) is admissible if and only if it factors as an admissible epic followed by an admissible monic.
\end{Lemma}

\begin{Proof}
  Obviously, the class of admissible morphisms is closed under composition, which proves one direction.

  Conversely, if $f\colon\mathbf{V}\to\mathbf{W}$ is admissible, it factors as $\mathbf{V}\twoheadrightarrow\im(f)\hookrightarrow\mathbf{W}$, where $\mathbf{V}\twoheadrightarrow\im(f)$ is an admissible epic ($\ker(f)\to\mathbf{V}\to\im(f)$ is a kernel-cokernel pair) and $\im(f)\hookrightarrow\mathbf{W}$ is an admissible monic ($\im(f)\hookrightarrow\mathbf{W}\to\coker(f)$ is a kernel-cokernel pair).
\end{Proof}

\begin{Lemma}\label{lem:inj-mon}
  Let $f\colon \mathbf{V}\to\mathbf{W}$ be an admissible morphism of F-zips over $k$.

  \begin{enumerate}
  \item The underlying map of sets of $f$ is injective if and only if $f$ is an admissible monic.
  \item The underlying map of sets of $f$ is surjective if and only if $f$ is an admissible epic.
  \end{enumerate}
\end{Lemma}

\begin{Proof}
  Essentially, this is true because it is true in ($k$-Mod) and the forgetful functor (F-zips$/k)\to(k$-Mod) preserves kernels and cokernels. In more detail:

  $f$ not injective $\implies \ker(f)\ne 0\implies f$ is not a monomorphism $\implies f$ is not a kernel.

  And if $f$ is injective, $\mathbf{V}\to\mathbf{W}\to\coker(f)$ is a kernel-cokernel pair (the isomorphism $V\cong\im(f)=\ker(W\to\coker(f))$ of vector spaces is an isomorphism of F-zips).

  Dually for surjectivity.
\end{Proof}

\begin{Proposition}
  (F-zips$/k$) is an exact category with the notion of admissible morphisms as in Definition \ref{def:adm}.
\end{Proposition}

\begin{Proof}
  We have shown that (F-zips$/k$) is additive in Lemma \ref{lem:addit}.

  Define $\mathscr{E}$ to be the class of kernel-cokernel pairs $(i,p)$ with $i$ and $p$ admissible in the sense of Definition \ref{def:adm}. Note that the definition of an admissible morphism relative to $\mathscr{E}$ (\cite{0811.1480}, Definition 8.1) then by Lemma \ref{lem:adm} will be equivalent to ours.

  Given a diagram
  \begin{center}
      \begin{tikzpicture}[node distance=4cm, auto]
        \node (grCV) at (0,2) {$\mathbf{V}$}; 
        \node (grCW) at (3,2) {$\mathbf{W}$};
        \node (grDW) at (3,0) {$\mathbf{W}'$};
        \node (grDV) at (0,0) {$\mathbf{V}'$};
    
        \draw[->] (grCV) to node {$f$} (grCW);
        \draw[->] (grDV) to node {$f'$} (grDW);
        \draw[->] (grCV) to node[swap] {$\cong$} (grDV);
        \draw[->] (grCW) to node {$\cong$} (grDW);
      \end{tikzpicture}
    \end{center}
    it is easy to see that $f$ is admissible if and only if $f'$ is, i.e., $\mathscr{E}$ is  closed under isomorphisms.
    
    We will now show that $\mathscr{E}$ satisfies the axioms of an exact structure.

    $\mathscr{E}$ obviously satisfies axioms [E0]--[E1$^\text{op}$] of an exact structure on (F-zips$/k$), that is, all identity morphisms are admissible monics (resp. admissible epics) and the class of admissible monics (resp. admissible epics) is closed under composition.

    For axioms [E2] and [E2$^\text{op}$] we will first need to demonstrate the existence  of the push-out of an admissible monic along an arbitrary morphism and, dually, the pull-back of an admissible epic along an arbitrary morphism.

    Let $i\colon\mathbf{V}\hookrightarrow\mathbf{W}$ be an admissible monic and let  $g\colon \mathbf{V}\to \mathbf{T}$ be any morphism. Then we claim that $z\colon\mathbf{V}\to\mathbf{W}\oplus\mathbf{T},\, v\mapsto(i(v),-g(v))$ is admissible. To see this, let $y\in z(V)\cap C^i$, $y=(i(v),-g(v))$, $v\in V$. Then $i(v)\in C^i$, hence $v\in C^i$ (by admissibility of $i$), which is what we had to show. The push-out of $i$ along $g$ then is $\coker(z)$. The resulting morphism $\mathbf{W}\to\coker(z)$ is the inclusion $\mathbf{W}\to \mathbf{W}\oplus\mathbf{T}$ followed by the cokernel of an admissible morphism, hence it is a composition of admissible morphisms, hence it is admissible. Also, push-out preserves exactness (i.e., the push-out of $i$ along $g$ is again monic), because of Lemma \ref{lem:inj-mon} and this being true in the category of $k$-vector spaces.

    Thus [E2] is proved. [E2$^\text{op}$] then is proved much in the same way.
\end{Proof}

\subsection{Exterior and symmetric powers}
\label{sec:aussere-und-symm}

We now wish to define the exterior powers of an F-zip. We already know how tensor products of F-zips work and of course we base  our construction on that. One obvious way to go about this might be to just divide out the same subspace that one divides out of the tensor product of $k$-vector spaces in order to obtain the exterior power of $k$-vector spaces. However, we would have no obvious F-zip structure on this quotient and for this reason we will instead construct the exterior powers as images of certain admissible morphisms in the category of F-zips.

But first we describe how this works for plain $R$-modules.

\begin{Definition} \label{ext-power}
  Let $R$ be a ring and $M$ an $R$-module.

  There is a natural action of $S_m$ on $M^{\otimes m}$: For $\pi\in S_m$ and $v_1,\dotsc,v_m\in M$ define
  \begin{equation*}
    \pi(v_1\otimes\dotsb\otimes v_m):=v_{\pi(1)}\otimes\dotsb\otimes v_{\pi(m)}.
  \end{equation*}
  and extend  linearly to obtain an automorphism $\pi$ of $M^{\otimes m}$.

  Define
  \begin{equation*}
    A_m\colon M^{\otimes m}\to M^{\otimes m},\; x\mapsto \sum_{\pi\in S_m}\sgn(\pi)\,\pi(x).
  \end{equation*}

  One might also interpret $A_m((-)^{\otimes m})$ as a functor $(R\text{-Mod})\to(R\text{-Mod})$, since ${f^{\otimes m}\colon M^{\otimes m}\to N^{\otimes m}}$ (where $f\colon M\to N$ is some $R$-linear map) restricts to an $R$-linear map ${A_m(M^{\otimes m})\to A_m(N^{\otimes m})}$.
\end{Definition}

\begin{Lemma}\label{Alternator}
  Let $R$ be a ring, $M$ a free $R$-module.

  Let $\{e_i\}_{i\in I}$ be a basis of $M.$ For $\underline i=(i_1,\dotsc,i_m)\in I^m$ write $e_{\underline i}:= e_{i_1}\otimes \dotsm \otimes e_{i_m}\in M^{\otimes m}.$

  \begin{enumerate}
  \item The kernel of $A_m$ is generated by all $v_1\otimes\dotsm\otimes v_m$, where $(v_1,\dotsc,v_m)\in M^m$ with $v_j=v_l$ for some $j\ne l$.
  \item $\bigwedge^m(M)\cong A_m(M^{\otimes m})$ naturally
  \item $A_m(M^{\otimes m})\subseteq\{\,z\in M^{\otimes m} \suchthat \forall \pi\in S_m\colon \pi(z)=\sgn(\pi)\,z\,\}$ (``alternating tensors are skew-symmetric'') and if $R\to R,\; x\to 2x$ is injective, equality holds.
  \item If $\sum_{\underline i\in I^m} a_{\underline i}e_{\underline i}\in A_m(M^{\otimes m})$, then $\sum_{\underline i\in I^m, i_1<\dotsb<i_m} a_{\underline i}e_{\underline i}$ is a preimage under $A_m.$
  \end{enumerate}
\end{Lemma}

\begin{Proof}
  Set $P:=\{\underline i\in I^m \suchthat i_j\ne i_l\text{ whenever } j\ne l\}.$

  1.+2. Let $(v_1,\dotsc,v_m)\in M^m$ with $v_j=v_l$ for some $j\ne l$ and $v:=v_1\otimes\dotsm\otimes v_m.$

  We have a partition
  $S_m= Z_< \sqcup Z_>$, where $Z_<=\{\pi\suchthat \pi(j) < \pi(l)\}$,
  ${Z_>=\{\pi \suchthat \pi(j)> \pi(l)\}}$, and
  \begin{equation*}
    \sum_{\pi\in Z_<}\sgn(\pi)\,\pi(v)=\sum_{\pi'\in Z_>}\sgn(\pi'\circ(j\quad l))\,(\pi'\circ(j\quad l))(v)=-\sum_{\pi'\in Z_>}\sgn(\pi')\,\pi'(v),
  \end{equation*}
  therefore $A_m(v)=0.$

  This shows that $N:=\langle v_1\otimes\dotsm\otimes v_m \suchthat \exists j\ne l\colon v_j=v_l\rangle\subseteq \ker(A_m).$

  For the converse let $x=\sum_{\underline i\in I^m}a_{\underline i}e_{\underline i}\in \ker(A_m).$ Since we already know that ${\sum_{\underline i\in I^m\setminus P}a_{\underline i}e_{\underline i}\in N}$, we may  assume that $a_{\underline i}=0$ for $\underline i\notin P$. 

Fix $\underline i\in P.$ $A_m(x)=0$ implies $\sum_{\pi\in S_m}\sgn(\pi)\,a_{\pi(\underline i)}=0.$

This is a linear equation, whose solution space is generated by all $f_1-\sgn(\pi)\,f_\pi$, $\pi\in S_m,$ where $\{f_\pi\}_\pi$ is the canonical basis of $R^{(S_m)}.$

Writing $(a_{\pi(\underline i)})_{\pi\in S_m}=\sum r_{\underline i,\pi} (f_1-\sgn(\pi)\, f_\pi),$ we obtain
\begin{equation*}
  x = \sum_{i_1<\dotsb < i_m} \sum_{\pi\in S_m}a_{\pi(\underline i)}e_{\pi(\underline i)} = \sum_{i_1<\dotsb < i_m} \sum_{\pi\in S_m}r_{\underline i,\pi}(e_{\underline i}-\sgn(\pi)\,e_{\pi(\underline i)}).
\end{equation*}

And $e_{\underline i} - \sgn(\pi)\,e_{\pi(\underline i)}$ clearly vanishes in $\bigwedge^mM=M^{\otimes m}/N.$

Hence $x\in N$ and $\ker(A_m)=N.$ By the the isomorphism theorem, $A_m$ gives rise to an isomorphism
\begin{equation*}
\bigwedge^m M=M^{\otimes m}/N\cong A_m(M^{\otimes m}).
\end{equation*}

(Note also that 4. will enable us to write down the inverse map explicitely.)

Naturality is trivial.

3. The inclusion ``$\subseteq$'' is trivial.

 Suppose that multiplication by $2$ is injective (in $R$ and thus in every free $R$-module).

  Let $z\in M^{\otimes m}$ with $\forall \pi\in
  S_m\colon\pi(z)=\sgn(\pi)\,z.$ Express $z$ as a linear combination of
  basis vectors $z =\sum_{\underline{i}} a_{\underline{i}} e_{\underline{i}}.$

  We have $a_{\pi^{-1}(\underline{i})}=\sgn(\pi)\,a_{\underline{i}}$ for all $\underline i,\pi$ by the
  assumption. So, if $i_j=i_l$, $j\ne l$ and $\tau=(j\quad l)\in S_m$, we have $a_{\underline i}=a_{\tau(\underline i)}=\sgn(\tau)\,a_{\underline i}=-a_{\underline i}$, that is, $2a_{\underline i}=0$, which means that $a_{\underline i}=0.$
  Hence
\begin{align*}
    z&=\sum_{\underline i\in P}a_{\underline i}e_{\underline i}
    =\sum_{i_1<\dotsb<i_m}\sum_{\pi\in
      S_m}a_{\pi(\underline i)}e_{\pi(\underline i)}=\sum_{i_1<\dotsb<i_m}\sum_{\pi\in
      S_m}\sgn(\pi^{-1})\,a_{\underline i}e_{\pi(\underline i)}\\
    &=\sum_{i_1<\dotsb<i_m}a_{\underline i}\sum_{\pi\in
      S_m}\sgn(\pi)\,e_{\pi(\underline i)}=A_m(\sum_{i_1<\dotsb<i_m}a_{\underline i}e_{\underline i})\in A_m(M^{\otimes m}).
  \end{align*}

  This calculation also proves 4.: If $z\in A_m(M^{\otimes m})$, we know by 1. that there is a preimage $x$ of $z$ under $A_m$ that can be written in terms of basis vectors $e_{\underline i}$ with $\underline i\in P.$ But then $z$ also can be written in terms of these vectors by definition of $A_m.$ Now carry out the calculation above.
\end{Proof}

\begin{Remark}
  The statements 1.--3. of the preceding lemma also hold if $M$ is only assumed to be \textit{locally} free (since the questions of whether one submodule is included in another and whether a natural isomorphism exists, are local in nature).
\end{Remark}

  Back to F-zips. We still want to construct exterior powers of F-zips and our course of action -- motivated by  statement 2. of the preceding lemma -- will be to show that $A_m$ is an admissible morphism of F-zips and then to define the $\bigwedge^m\mathbf{V}$ as the image of $A_m\colon \mathbf{V}^{\otimes m}\to \mathbf{V}^{\otimes m}$.

\begin{Lemma}\label{ext-power-adm}
 $A_m\colon V^{\otimes m}\to V^{\otimes m}$ is an admissible morphism of F-zips.
\end{Lemma}

\begin{Proof}
 We have
  \begin{align*}
    A_m(C^r(\mathbf{V}^{\otimes m})) &= A_m(\sum_{i_1+\dotsb+i_m=r}\bigotimes_{\nu=1}^mC^{i_\nu}\mathbf{V})\subseteq \sum_{i_1+\dotsb+i_m=r}\bigotimes_{\nu=1}^mC^{i_\nu}\mathbf{V} = C^r(\mathbf{V}^{\otimes m}), \\
    A_m(D_r(\mathbf{V}^{\otimes m})) &= A_m(\sum_{i_1+\dotsb+i_m=r}\bigotimes_{\nu=1}^mD_{i_\nu}\mathbf{V})\subseteq \sum_{i_1+\dotsb+i_m=r}\bigotimes_{\nu=1}^mD_{i_\nu}\mathbf{V} = D_r(\mathbf{V}^{\otimes m}).
  \end{align*}

  Next, we show
  \begin{align*}
    \gr_r^DA_m\circ(\bigoplus_{i_1+\dotsb+i_m=r}\bigotimes_{\nu=1}^m\varphi_{i_\nu})
    &= (\bigoplus_{i_1+\dotsb+i_m=r}\bigotimes_{\nu=1}^m\varphi_{i_\nu}) \circ \gr_C^rA_m.
  \end{align*}

  For $c^{i}_j\in C^{i_j}\mathbf{V}$, we have
  \begin{align*}
    &\gr_C^rA_m(\sum_{i_1+\dotsb+i_m=r} c^{i}_1\otimes\dotsm\otimes c^{i}_m+C^{r+1}(\mathbf{V}^{\otimes m})) \\
    &= \sum_{i_1+\dotsb+i_m=r} A_m(c^{i}_1\otimes\dotsm\otimes c^{i}_m)+C^{r+1}(\mathbf{V}^{\otimes m}) \\
    &= \sum_{i_1+\dotsb+i_m=r} \sum_{\pi\in S_m}\sgn(\pi)\,c^{i}_{\pi(1)}\otimes\dotsm\otimes c^{i}_{\pi(m)}+C^{r+1}(\mathbf{V}^{\otimes m}) \\
    &\mathop{\widehat=} \left(\sum_{\pi\in S_m}\sgn(\pi)\,(c^{\pi^{-1}(i)}_{\pi(1)}+C^{i_1+1}\mathbf{V})\otimes\dotsm\otimes (c^{\pi^{-1}(i)}_{\pi(m)}+C^{i_m+1}\mathbf{V})\right)_{i_1+\dotsb+i_m=r},
  \end{align*}
  so that, choosing a representative $\varphi_\nu^\mathrm{Rep}(\bar c)\in D_\nu\mathbf{V}$ of $\varphi_\nu(\bar c)\in D_\nu\mathbf{V}/D_{\nu-1}\mathbf{V}$, we have
  \begin{align*}
    &((\bigoplus_{i_1+\dotsb+i_m=r}\bigotimes_{\nu=1}^m\varphi_{i_\nu}) \circ \gr_C^rA_m(\sum_{i_1+\dotsb+i_m=r} c^{i}_1\otimes\dotsm\otimes c^{i}_m+C^{r+1}(\mathbf{V}^{\otimes m})))\\
    &=\left(\sum_{\pi\in S_m}\sgn(\pi)\,\varphi_{i_1}(c^{\pi^{-1}(i)}_{\pi(1)}+C^{i_1+1}\mathbf{V})\otimes\dotsm\otimes \varphi_{i_m}(c^{\pi^{-1}(i)}_{\pi(m)}+C^{i_m+1}\mathbf{V})\right)_{i_1+\dotsb+i_m=r} \\
    &\mathop{\widehat=} \sum_{i_1+\dotsb+i_m=r}\sum_{\pi\in S_m}\sgn(\pi)\,\varphi_{i_1}^\mathrm{Rep}(c^{i}_{\pi(1)}+C^{i_1+1}\mathbf{V})\otimes\dotsm\otimes \varphi_{i_m}^\mathrm{Rep}(c^{i}_{\pi(m)}+C^{i_m+1}\mathbf{V})+D_{r-1}(\mathbf{V}^{\otimes m}) \\
    &= \sum_{i_1+\dotsb+i_m=r}A_m(\varphi_{i_1}^\mathrm{Rep}(c^{i}_{1}+C^{i_1+1}\mathbf{V})\otimes\dotsm\otimes \varphi_{i_m}^\mathrm{Rep}(c^{i}_{m}+C^{i_m+1}\mathbf{V}))+D_{r-1}(\mathbf{V}^{\otimes m}) \\
    &= \gr_r^DA_m(\sum_{i_1+\dotsb+i_m=r}\varphi_{i_1}^\mathrm{Rep}(c^{i}_{1}+C^{i_1+1}\mathbf{V})\otimes\dotsm\otimes \varphi_{i_m}^\mathrm{Rep}(c^{i}_{m}+C^{i_m+1}\mathbf{V}) + C^{r+1}(\mathbf{V}^{\otimes m})) \\
    &\mathop{\widehat=} \gr_r^DA_m(\left(\varphi_{i_1}(c^{i}_{1}+C^{i_1+1}\mathbf{V})\otimes\dotsm\otimes \varphi_{i_m}(c^{i}_{m}+C^{i_m+1}\mathbf{V})\right)_{i_1+\dotsb+i_m=r}) \\
    &= \gr_r^DA_m\circ(\bigoplus_{i_1+\dotsb+i_m=r}\bigotimes_{\nu=1}^m\varphi_{i_\nu})(\left((c^{i}_{1}+C^{i_1+1}\mathbf{V})\otimes\dotsm\otimes (c^{i}_{m}+C^{i_m+1}\mathbf{V})\right)_{i_1+\dotsb+i_m=r}) \\
    &=\gr_r^DA_m\circ(\bigoplus_{i_1+\dotsb+i_m=r}\bigotimes_{\nu=1}^m\varphi_{i_\nu})(\sum_{i_1+\dotsb+i_m=r} c^{i}_1\otimes\dotsm\otimes c^{i}_m+C^{r+1}(\mathbf{V}^{\otimes m})).
  \end{align*}
  
  We conclude that $A_m\in\Hom_{(\text{F-zips/}k)}(\mathbf{V}^{\otimes m},\mathbf{V}^{\otimes m})$ and still have to prove admissibility.

  Choose  a basis $\{e_i\}_i$ of $V$ adapted to the filtration $C^\bullet\mathbf{V}.$
  Suppose that ${z\in C^r(\mathbf{V}^{\otimes m})\cap A_m(V^{\otimes m})}$. Then we
  may write $z=\sum_{\underline i} a_{\underline i} e_{\underline i}$
  with $e_{i_j}\in C^{s(\underline i,j)}(\mathbf{V}),\;
  \sum_js(\underline i,j)=r,$ and, according to Lemma \ref{Alternator}, $\sum_{i_1
    < \dotsb < i_m}a_{\underline i}e_{\underline i}$ is a preimage
  under $A_m.$ It follows that $A_m(\mathbf{V}^{\otimes m})\cap
  C^r(\mathbf{V}^{\otimes m})\subseteq A_m(C^r(\mathbf{V}^{\otimes
    m}))$ and, similarly, $A_m(\mathbf{V}^{\otimes m})\cap
  D_r(\mathbf{V}^{\otimes m})\subseteq A_m(D_r(\mathbf{V}^{\otimes
    m}))$, i.e., $A_m$ is admissible.  
\end{Proof}

\begin{Definition}
  Define
  \begin{align*}
    \bigwedge^m\mathbf{V} &:=\im A_m\\
    &= (\bigwedge^mV, \bigwedge^mV\cap C^\bullet(\mathbf{V}^{\otimes m}),
    \bigwedge^mV\cap D_\bullet(\mathbf{V}^{\otimes m}),
    \left.\varphi_\bullet(\mathbf{V}^{\otimes
        m})\right|_{\bigwedge^m\mathbf{V}}).
  \end{align*}
\end{Definition}

We give a short account of some basic facts and identities concerning the exterior powers, which will be useful later on.

\begin{Remark} \label{Cbigwedge}
  One has
  \begin{align*}
  C^r(\bigwedge^m\mathbf{V})&=\langle v_1\wedge\dotsm\wedge v_m \suchthat v_j\in C^{i_j}\mathbf{V},\; i_1+\dotsb+i_m=r\rangle \\
  &=\langle v_1\wedge\dotsm\wedge v_m \suchthat v_j\in C^{i_j}\mathbf{V},\; i_1+\dotsb+i_m=r,\; i_1\geq \dotsb\geq i_m\rangle
\end{align*}
\end{Remark}

\begin{Lemma}
  Let $A$ be a ring, $\left(M_\lambda\right)_{\lambda\in L}$ a family of $A$-modules.

  Then
  \begin{equation*}
  \bigwedge(\bigoplus_{\lambda \in L} M_\lambda) \cong  {^g}\bigotimes_{\lambda \in L} \bigwedge M_\lambda
  \end{equation*}
  as graded $A$-algebras.
\end{Lemma}

  Here, if $B_\lambda=\bigoplus_i B_\lambda^i$, $\lambda\in L$, are graded $A$-algebras, then ${^g}\bigotimes_{\lambda \in L} B_\lambda$ is the skew tensor  algebra of $\left(B_\lambda\right)_\lambda$ that is characterized by the following universal property: There are morphisms $f_\lambda\colon B_\lambda\to {^g}\bigotimes_{\lambda\in L} B_\lambda$ such that if $T$ is a graded  $A$-algebra and $g_\lambda\colon B_\lambda\to T$ are  morphisms satisfying $g_\lambda(x)\cdot g_\mu(y)=(-1)^{kl}g_\mu(y)\cdot g_\lambda(x)$ for all $x\in B_\lambda^k$ and $y\in B_\mu^l$, then there exists a unique $g\colon {^g}\bigotimes_{\lambda\in L} B_\lambda\to T$ with $g\circ f_\lambda=g_\lambda$ (cf. \cite{bourbaki1998algebra}, III, §4.7, prop. 10).

  The $r$-th graded piece of ${^g}\bigotimes_{\lambda \in L} B_\lambda$ is given by $\bigoplus B_\lambda^{i_\lambda}$, where the direct sum is taken over all $\left(i_\lambda\right)_{\lambda\in L}$ with $\sum i_\lambda=r$, and the morphisms $f_\lambda$ from above are the obvious ones.

\begin{Proof}
  (\cite{bourbaki1998algebra}, Algebra, III, §7.7, prop. 10)

  The natural maps $j_\lambda\colon M_\lambda\to \bigoplus_{\lambda \in L} M_\lambda$ give rise to
  \begin{equation*}
    \bigwedge j_\lambda\colon\bigwedge M_\lambda \to \bigwedge \bigoplus_{\lambda \in L} M_\lambda,
  \end{equation*}
  using functoriality of the exterior algebra.

  We have $(\bigwedge j_\lambda)(x)\wedge(\bigwedge j_\mu)(y) = (-1)^{kl}(\bigwedge j_\mu)(y)\wedge(\bigwedge j_\lambda)(x)$ for homogeneous elements $(x,y) \in \bigwedge^k M_\lambda \times \bigwedge^l M_\mu$. By the universal property of the skew tensor algebra, we get a homomorphism
  \begin{equation*}
    g\colon {^g}\bigotimes_{\lambda\in L}\bigwedge M_\lambda \to \bigwedge \bigoplus_{\lambda \in L} M_\lambda \quad \text{ with } \quad \bigwedge j_\lambda = g\circ f_\lambda \quad \forall\lambda\in L,
  \end{equation*}
  where $f_\lambda\colon \bigwedge M_\lambda \to {^g}\bigotimes_{\lambda\in L}\bigwedge M_\lambda$ is the natural homomorphism.

This is the desired isomorphism.
\end{Proof}

\begin{Corollary} \label{ext-coprod}
  Let $A$ be a ring, $\left(M_\lambda\right)_{\lambda\in L}$ a family of $A$-modules.
  \begin{equation*}
  \bigwedge^m(\bigoplus_{\lambda \in L} M_\lambda) \cong  \bigoplus_{\sum i_\lambda=m}\bigotimes_{\lambda \in L} \bigwedge^{i_\lambda} M_\lambda
  \end{equation*}
  as $A$-modules.
\end{Corollary}

\begin{Remark}
  If $L=\{1,2\}$, the isomorphism from Corollary \ref{ext-coprod} is given by
  \begin{align*}
    (v_1,0)\wedge\dotsm\wedge(v_i,0)\wedge(0,w_{i+1})\wedge\dotsm\wedge(0,w_m)&\mapsto (v_1\wedge\dotsm\wedge v_i)\otimes (w_{i+1}\wedge\dotsm\wedge w_m)
  \end{align*}
  and $\bigwedge^m(M_1\oplus M_2)$ is generated by vectors of the form ${(v_1,0)\wedge\dotsm\wedge(v_i,0)\wedge(0,w_{i+1})\wedge\dotsm\wedge(0,w_m)}$.
\end{Remark}

\begin{Remark}\label{rem:iso-of-f-zips}
  Let $\Phi\colon \bigwedge^m(V\oplus W)\xrightarrow{\sim} \bigoplus_{a\geq 0}\bigwedge^aV\otimes\bigwedge^{m-a}W$ be the isomorphism from Corollary \ref{ext-coprod}.

  Then $\Phi$ is an isomorphism of F-zips,
  \begin{equation*}
    \bigwedge^m(\mathbf{V}\oplus\mathbf{W})\cong\bigoplus_{i+j=m}\bigwedge^i\mathbf{V}\otimes\bigwedge^j\mathbf{W}.
  \end{equation*}
\end{Remark}

\begin{Proof}
  \begin{align*}
    C^r\bigwedge^m(\mathbf{V}\oplus\mathbf{W})
    &=\left\langle\bigwedge_{j=1}^m(v_j,w_j)\suchthat v_j\in C^{i_j}\mathbf{V},w_j\in C^{i_j}\mathbf{W},\;\sum_{j=1}^m i_j=r\right\rangle \\
    &=\left\langle\bigwedge_{j=1}^m(v_j,w_j)\suchthat \forall j\colon(v_j\in C^{i_j}\mathbf{V},w_j\in C^{i_j}\mathbf{W},\; v_j=0 \text{ oder } w_j=0),\;\sum_{j=1}^m i_j=r\right\rangle \\
    &=\sum_{a\geq 0}\left\langle\bigwedge_{j=1}^a(v_j,0)\wedge\bigwedge_{j=a+1}^m(0,w_j)\suchthat v_j\in C^{i_j}\mathbf{V},\;w_j\in C^{i_j}\mathbf{W},\;\sum_{j=1}^m i_j=r\right\rangle,
  \end{align*}
  hence
  \begin{align*}
    \Phi(C^r\bigwedge^m(\mathbf{V}\oplus\mathbf{W}))
    &=\bigoplus_{a\geq 0}\left\langle\bigwedge_{j=1}^av_j\otimes\bigwedge_{j=a+1}^mw_j\suchthat v_j\in C^{i_j}\mathbf{V},\;w_j\in C^{i_j}\mathbf{W},\;\sum_{j=1}^m i_j=r\right\rangle\\
    &=\bigoplus_{a\geq 0}\sum_{s+t=r}\left\langle\bigwedge_{j=1}^av_j\otimes\bigwedge_{j=a+1}^mw_j\suchthat v_j\in C^{i_j}\mathbf{V},\;w_j\in C^{i_j}\mathbf{W},\;\sum_{j=1}^a i_j=s,\;\sum_{j=a+1}^m i_j=t\right\rangle\\
    &=\bigoplus_{a\geq 0}\sum_{s+t=r}\left\langle\bigwedge_{j=1}^av_j\otimes\bigwedge_{j=a+1}^mw_j\suchthat \bigwedge_{j=1}^av_j\in C^s\bigwedge^a\mathbf{V},\; \bigwedge_{j=a+1}^mw_j\in C^t\bigwedge^{m-a}\mathbf{W}\right\rangle\\
    &=\bigoplus_{a\geq 0}\sum_{s+t=r}C^s\bigwedge^a\mathbf{V}\otimes C^t\bigwedge^{m-a}\mathbf{W}\\
    &=\bigoplus_{a\geq 0}C^r(\bigwedge^a\mathbf{V}\otimes\bigwedge^{m-a}\mathbf{W})\\
    &=C^r(\bigoplus_{a\geq 0}\bigwedge^a\mathbf{V}\otimes\bigwedge^{m-a}\mathbf{W})
  \end{align*}
  and simlarly for $D_\bullet$.

  Also,
  \begin{equation*}
    \left(\bigoplus_{a\geq 0}\bigoplus_{i+j=r}\varphi_i(\bigwedge^a\mathbf{V})\otimes\varphi_j(\bigwedge^{m-a}\mathbf{W})\right)\circ\gr_C^r\Phi=\gr_r^D\Phi\circ \varphi_r(\bigwedge^m(\mathbf{V}\oplus\mathbf{W})),
  \end{equation*}
  as a cumbersome calculation similar to the one in \ref{ext-power-adm} shows.
\end{Proof}

Finally, we note that we also have the similar notion of \textit{symmetric} powers of F-zips.

\begin{Definition}
  Define
  \begin{equation*}
    B_m\colon \bigoplus_{\pi\in S_m} V^{\otimes m} \to V^{\otimes m},\;
    \left(x_\pi\right)_\pi \mapsto \sum_{\pi\in S_m}(\pi(x_\pi)-x_\pi).
  \end{equation*}

  We then have
  \begin{align*}
    \im B_m &=\left\langle\pi(x)-x\suchthat x\in V^{\otimes m},\;\pi\in S_m\right\rangle \\
    &=\left\langle\tau(x)-x\suchthat x\in V^{\otimes m},\;\tau\in S_m \text{ transposition of two adjacent elements}\right\rangle.
  \end{align*}

  $B_m$ is an admissible morphism of F-zips and
  \begin{equation*}
    S^m\mathbf{V}:=\coker(B_m).
  \end{equation*}
\end{Definition}

\subsection{Base change}
\label{sec:korperwechsel}

\begin{Definition}
  Let $k'/k$ be a field extension of perfect fields. We won't distinguish between the Frobenius endomorphism of $k$ and the Frobenius endomorphism of $k'$.

  Then, defining $a(v\otimes b):= v\otimes (ab)$, $V\otimes_kk'$ is a $k'$-vector space  and we  have ${(C^i\otimes_kk')/(C^{i+1}\otimes_kk')\cong (C^i/C^{i+1})\otimes_kk'}$.

  Define
  \begin{equation*}
    \mathbf{V}_{k'}:=(V\otimes_kk', C^\bullet\otimes_kk', D_\bullet\otimes_kk', \varphi_\bullet\otimes\sigma).
  \end{equation*}
  This (evidently) is an F-zip over $k'$.
\end{Definition}

\begin{Lemma}
  $\bigwedge^m(\mathbf{V}_{k'})\cong\left(\bigwedge^m\mathbf{V}\right)_{k'}$.
\end{Lemma}

\begin{Proof}
  As in Remark \ref{rem:iso-of-f-zips}: The natural isomorphism of vector spaces respects the F-zip structure.
\end{Proof}

\subsection{Dual}
\label{sec:dual}
\begin{Definition}
  Define $\mathbf{V}^\vee$ by
  \begin{align*}
    C^i(\mathbf{V}^\vee) &:=  (V/C^{1-i}(\mathbf{V}))^\vee=\{\lambda\in V^\vee\suchthat C^{1-i}(\mathbf{V})\subseteq\ker \lambda\}, \\
    D_i(\mathbf{V}^\vee) &:=  (V/D_{-1-i}(\mathbf{V}))^\vee=\{\lambda\in V^\vee\suchthat D_{-1-i}(\mathbf{V})\subseteq\ker \lambda\}, \\
    \varphi_i(\mathbf{V}^\vee)&:=((\varphi_{-i}(\mathbf{V}))^{-1})^\vee.
  \end{align*}

  We have
  \begin{align*}
    \gr_C^i(\mathbf{V}^\vee)&\cong\gr_C^{-i}(\mathbf{V})^\vee, \\
    \gr_i^D(\mathbf{V}^\vee)&\cong\gr_{-i}^D(\mathbf{V})^\vee.
  \end{align*}
\end{Definition}

\section{Classification of F-zips}
\label{sec:klassifikation-von-f-zips}

In this section we will give a system of representatives for the isomorphism classes of F-zips of a fixed type over a fixed algebraically closed field.

Fix a type $\tau\colon \Z\to\N_0$. Let $\supp(\tau)=\{i_1>\dotsb>i_r\}$ and $n_\nu:=\tau(i_\nu)$, $n:=\sum_{\nu=1}^rn_\nu$ and $m_\mu:=\sum_{\nu=1}^\mu n_\nu$. Set
\begin{align*}
  W &:= S_n, \\
  I &:= \{(i\quad i+1)\in S_n\suchthat i\in\{1,\dotsc,n-1\}\}, \\
  J &:= \{(i\quad i+1)\in  I\suchthat i\notin\{m_\mu\}_{\mu=1,\dotsc,r-1}\}, \\
  W_J &:= (\text{the subgroup of } W \text{ generated by } J)=\prod_{\mu=1}^{r}S(\{m_{\mu-1}+1,\dotsc,m_{\mu}\})\subseteq S_n,
\end{align*}
where the inclusion is given by
\begin{equation*}
  \prod_{\mu=1}^{r}S(\{m_{\mu-1}+1,\dotsc,m_{\mu}\})\ni(\pi_1,\dotsc,\pi_{r})\mapsto \pi_1\dotsb\pi_{r}\in S(\{1,\dotsc,n\})=S_n
\end{equation*}
($Y\subseteq X \implies S(Y)\subseteq S(X)$ in the obvious way).

\begin{Definition}[Length of a permutation] Let $w\in S_n$.
  \begin{equation*}
    \ell(w):=\#(\text{inversions of } w)=\#\{(i,j)\in\{1,\dotsc,n\}^2\suchthat i<j,\; w(i)>w(j)\}
  \end{equation*}
\end{Definition}

Let $^JW$ be the set of the unique $\ell$-minimal representatives of the left cosets  $W_J\verb|\|W=\{W_Jw\suchthat w\in W\}$. We will now describe $^JW$ and in particular justify the usage of the word ``unique'' in the preceding sentence.

Fix $w\in W$. Representatives of $W_Jw$ are of the form
\begin{equation*}
  \tilde w=w_1\dotsm w_{r}w, \quad\text{where } w_\mu\in S(\{m_{\mu-1}+1,\dotsc,m_{\mu}\})\subseteq S_n.
\end{equation*}

Consequently,  inversions $i<j,\;w(i)>w(j)$ with $w(i),w(j)\in\{m_{\mu-1}+1,\dotsc,m_{\mu}\}$ for some $\mu\in\{1,\dotsc,r\}$ are ``fixable'' and other inversions are not.

Thus, the $\ell$-minimal representative $\tilde w$ of $W_Jw$ satisfies
\begin{equation}
  \label{eq:min-rep}
  \tilde w^{-1}(m_{\mu-1}+1)<\dotsb<\tilde w^{-1}(m_{\mu})\qquad \forall \mu=1,\dotsc,r
\end{equation}
and
\begin{equation*}
  ^JW=\{\tilde w\in W\suchthat \eqref{eq:min-rep}\text{ holds for }\tilde w\}.
\end{equation*}

Set
\begin{align*}
  w_0(i)&:=n+1-i, \\
  w_{0,J}(i)&:=m_{\mu}+m_{\mu-1}+1-i\quad (m_{\mu-1}+1\leq i\leq m_{\mu}), \\
  w_0^J&:=w_0w_{0,J}, \quad w_0^J(i)=n-m_{\mu}-m_{\mu-1}+i \quad (m_{\mu-1}+1\leq i\leq m_{\mu}), \\
  w^\circ&:=ww_0^J
\end{align*}

Then $w_0$ is the single longest permutation in $W$ and $w_{0,J}$ is the single longest permutation in $W_J$.

Define an F-zip $^w\mathbf{V}=(\F_p^n, C^\bullet, D_\bullet, \varphi_\bullet)$ of type $\tau$ over $\F_p$ as follows:

$\F_p^n=\langle e_1,\dotsc,e_n\rangle_{\F_p}$

Because of the type that we fixed, we necessarily need to have $C^i/C^{i+1}=D_i/D_{i-1}=0$   for $i\notin\{i_\nu\}_\nu$, so that $C^i=C^{i+1}$ and $D_i=D_{i-1}$ for $i\notin\{i_\nu\}_\nu$.  Hence, it is sufficient to define $C^{i_\nu}$ and $D_{i_\nu}$ for $\nu\in\{1,\dotsc,r\}$.

Set
\begin{align*}
  C^{i_\nu}&:=\langle e_1,\dotsc,e_{m_\nu}\rangle_{\F_p}, \\
  D_{i_\nu}&:=\langle e_{w(m_{\nu-1}+1)},\dotsc,e_{w(n)}\rangle_{\F_p}.
\end{align*}

Then
\begin{align*}
  C^{i_\nu}/C^{i_{\nu}+1}&=C^{i_\nu}/C^{i_{\nu-1}}=\langle e_{m_{\nu-1}+1},\dotsc,e_{m_\nu}\rangle_{\F_p}, \\
  D_{i_\nu}/D_{i_{\nu}-1}&=D_{i_\nu}/D_{i_{\nu+1}}=\langle e_{w(m_{\nu-1}+1)},\dotsc,e_{w(m_\nu)}\rangle_{\F_p}.
\end{align*}

Define $\varphi_\nu$ by $\varphi_\nu(e_l):=e_{w(l)}$.

\begin{Proposition} The map
  \begin{align*}
    ^JW&\to \{\text{F-zips/}k \text{ of type } \tau\}/{\cong},\\
    w &\mapsto [\left(^{w^\circ}\mathbf{V}\right)_k]_{\cong}
  \end{align*}
  is bijective, if $k$ is algebraically closed.
\end{Proposition}

\begin{Proof}
  \cite{goettingen}, Theorem 3.6
\end{Proof}

\subsection{The example of full F-zips}

\begin{Lemma}
Let $k$ be a field and let $V$ be an $n$-dimensional $k$-vector space. Let $V=F_n\supsetneq F_{n-1} \supsetneq \dotsb \supsetneq F_0 = 0$ and $V = G_n \supsetneq G_{n-1} \supsetneq \dotsb \supsetneq G_0 = 0$ be two full flags on $V$.

Then there exists a basis $\{e_i\}_{1\leq i\leq n}$ of $V$ and a permutation $w\in S_n$ with ${F_j= \langle e_i \suchthat i\leq j\rangle}$ and $G_j=\langle e_{w(i)} \suchthat i\leq j\rangle$ for all $1\leq i\leq n.$
\end{Lemma}

\begin{Proof}
First choose some basis $\{e_i\}_{1\leq i\leq n}$ of $V$ with $F_j= \langle e_i \suchthat i\leq j\rangle$ for all $j$.

Then $G_1$ is generated by $\sum \lambda_ie_i$ for some $\lambda_i\in k$. Let $w(1):=\ell$ be the largest index with $\lambda_\ell\ne 0$ and replace $e_\ell$ by $\sum \lambda_ie_i$. Then $\{e_i\}_{1\leq i\leq n}$ still is a basis of $V$ with $F_j= \langle e_i \suchthat i\leq j\rangle$ for all $j$.

Now $G_2$ is generated by $e_{w(1)}$ and $\sum \lambda_ie_i$ for some $\lambda_i\in k$, $\lambda_{w(1)}=0$.  Let $w(2):=\ell$ be the largest index with $\lambda_\ell\ne 0$, replace $e_\ell$ by $\sum \lambda_ie_i$ and proceed like this.
\end{Proof}

\begin{Proposition}
Let $k$ be an algebraically closed field.

If $\forall i\in\Z\colon \tau(i)\leq 1,$ then an F-zip $\mathbf{V}$ over $k$ of type $\tau$ is already determined by $(V,C^\bullet,D_\bullet)$ up to isomorphism.
\end{Proposition}

\begin{Proof}
Using the previous remark we choose a basis $\{e_i\}_{1\leq i\leq n}$ of  $(V,C^\bullet,D_\bullet).$ Let ${w\in S_n}$ be the corresponding permutation.

Let $\psi_\bullet$ be the family of isomorphism $\gr_C^r\to\gr_r^D$ with $e_i\mapsto e_{w(i)}.$

The family $\varphi_\bullet$ associated to $\mathbf{V}$ is given $e_i\mapsto \lambda_i e_{w(i)}$ for some $\lambda_i\in k^\times.$

We shall now determine a tuple $\left(a_i\right)_{1\leq i\leq n}\in (k^\times)^n$ such that the $k$-linear map ${f\colon V\to V,\; e_i\mapsto a_ie_i}$ yields an isomorphism between $(V,C^\bullet,D_\bullet,\psi_\bullet)$ and $\mathbf{V}=(V,C^\bullet,D_\bullet,\varphi_\bullet).$

For this to be true, the only conditions on $\{a_i\}_i$ (other than $a_i\in k^\times$ for all $i$) are that for every $i\in \Z$ the diagram
  \begin{center}
      \begin{tikzpicture}[node distance=4cm, auto]
        \node (grC1) at (0,2) {$\langle e_i\rangle$}; 
        \node (grCV) at (6,2) {$\langle e_i\rangle$};
        \node (grDV) at (6,0) {$\langle e_{w(i)}\rangle$};
        \node (grD1) at (0,0) {$\langle e_{w(i)}\rangle$};
    
        \draw[->] (grC1) to node {$f$} (grCV);
        \draw[->] (grD1) to node {$f$} (grDV);
        \draw[->] (grC1) to node {$e_i\mapsto e_{w(i)}$} node[swap] {$\sigma$-linear} (grD1);
        \draw[->] (grCV) to node[swap] {$e_i\mapsto \lambda_ie_{w(i)}$} node {$\sigma$-linear} (grDV);
      \end{tikzpicture}
    \end{center}
commutes, that is, $\forall i\in \Z\colon a_i^p\lambda_i=a_{w(i)}.$

Write $w=c_1\dotsm c_s$ as a product of disjoint cycles, where $c_i=(j_{i,1}\quad \dotsb \quad j_{i,m_i})$ ($m_i=\ord(c_i)$), such that $\{1,\dotsc,n\}$ is the set of all $j_{i,l}.$

Our conditions then become
\begin{align*}\tag{$\ast$}\label{cyccond}
	\lambda_{j_{i,1}}a_{j_{i,1}}^p &= a_{w(j_{i,1})}, \\
	\lambda_{w(i)}a_{w(j_{i,1})}^p &= a_{w^2(j_{i,1})}, \\
	&\dotsc, \\
	\lambda_{w^{m_i-1}(j_{i,1})}a_{w^{m_i-1}(j_{i,1})}^p &= a_{j_{i,1}}.
\end{align*}

So, choosing a non-trivial (i.e., invertible) root $a_{j_{i,1}}$ of
\begin{align*}
  X - \lambda_{w^{m_i-1}(j_{i,1})}\left(\lambda_{w^{m_i-2}(j_{i,1})}\left(\dotsb\left(\lambda_{w(j_{i,1})}\left(\lambda_{j_{i,1}}X^p\right)^p\right)^p \dotsb \right)^p\right)^p \\
  = X\left(1-\left(\prod_{l=1}^{m_i} \lambda_{w^{m_i-l}(j_{i,1})}^{p^{l-1}}\right)X^{p^{m_i}-1}\right)
\end{align*}
for all $i$ and defining $a_{w^l(j_{i,1})}$ successively by \eqref{cyccond}, will do it.
\end{Proof}

\begin{Definition}
  An F-zip as in the previous proposition is called a \textit{full} F-zip.
\end{Definition}

\section{The permutation associated to $\bigwedge^m(^w\mathbf{V})_k$}
\label{sec:perm-assoc-}

Let $k$ be algebraically closed and consider the F-zip $(^w\mathbf{V})_k$ of type $\tau$ associated to the permutation $w\in S_n$, $n=\sum_{i\in\Z} \tau(i)$. We now wish to determine $w'\in S_{n'}$ ($n'=\binom{n}{m}$) such that $\bigwedge^m(^w\mathbf{V})_k$ is associated to the permutation $w'\in S_{n'}$.

Set
\begin{equation*}
  T_m:=\{S\subseteq \{1,\dotsc,n\} : |S|=m\}
\end{equation*}
and for $\mu_\nu\in\{1,\dotsc,n\}$ set
\begin{equation*}
  E_{\{\mu_1,\dotsc,\mu_m\}}:=e_{\mu_1}\wedge\dotsm\wedge e_{\mu_m}\in \bigwedge^mk^n,
\end{equation*}
which of course is only well defined up to sign. For our purposes we can mostly ignore this inconvenience since we have the following lemma:

\begin{Lemma}\label{lem:vorzeichen-egal}
  Let $\mathbf{V},\mathbf{W}$ be two F-zips (over an algebraically closed field $k$) of arbitrary type $\tau$, whose underlying vector spaces are identical and have a basis $e_1,\dotsc,e_n$ adapted to both the ascending and the descending filtration, such that the isomorphisms $\gr_C^r\to\gr_r^D$ are given by $e_i\mapsto \mu_ie_{w(i)}$ and $e_i\mapsto \lambda_ie_{w(i)}$ respectively for some $w\in S_n$ and $\mu_i,\lambda_i\in k^\times.$ Then $\mathbf{V}\cong\mathbf{W}.$
\end{Lemma}

\begin{Proof}
  Without loss of generality, $\forall i\colon \mu_i=1$.

  An isomorphism is given by $f\colon V\to V,\; e_i\mapsto a_ie_i$ for appropriate $a_i\in k^\times,$ which we may determine as follows:

  We want the diagram
\begin{center}
      \begin{tikzpicture}[node distance=4cm, auto]
        \node (grCrV) at (0,2) {$\gr_C^r\mathbf{V}$}; 
        \node (grCrW) at (4,2) {$\gr_C^r\mathbf{W}$};
        \node (grDV) at (0,0) {$\gr_r^D\mathbf{V}$};
        \node (grDW) at (4,0) {$\gr_r^D\mathbf{W}$};
    
        \draw[->] (grCrV) to node {$f$} (grCrW);
        \draw[->] (grDV) to node {$f$} (grDW);
        \draw[->] (grCrW) to node {$e_i\mapsto \lambda_i e_{w(i)}$} (grDW);
        \draw[->] (grCrV) to node[swap] {$e_i\mapsto e_{w(i)}$} (grDV);
      \end{tikzpicture}
    \end{center}
    to commute, where $\gr_C^r\mathbf{V}=\gr_C^r\mathbf{W}=\langle e_i\suchthat i\in I_r\rangle$ for some $I_r\subseteq\{1,\dotsc,n\}$ and horizontal arrows are $k$-linear, vertical arrows are $\sigma$-linear.

    So, for  $i\in I_r$ we want the following equation to be true:
    \begin{align*}
      a_i^p\lambda_ie_{w(i)}&=\left(\left(e_i\mapsto\lambda_ie_{w(i)}\right)\circ f\right)(e_i) \\
      &\stackrel{!}= \left(f\circ\left(e_i\mapsto e_{w(i)}\right)\right)(e_i)=a_ie_{w(i)},
    \end{align*}
    that is, we need to choose $a_i\in k^\times$ such that $a_i^{p-1}=\lambda_i^{-1}.$ Since $k$ is algebraically closed, this is possible.
\end{Proof}

Then by definition, with respect to the basis $\left(E_S\right)_{S\in T_m}$ (which is adapted to both filtrations), $\bigwedge^m(^w\mathbf{V})_k$ is given by the permutation $\{\mu_\nu\}_\nu\mapsto \{w(\mu_\nu)\}_\nu$, that is, by $\iota(w)$ where $\iota\colon S_n\to S(T_m)$ is the obvious group action of $S_n$ on $T_m.$

Put differently,
\begin{equation*}
  \bigwedge^m(^w\mathbf{V})_k=(^{\iota(w)}\mathbf{V})_k.
\end{equation*}

Of course, $\bigwedge^m(^w\mathbf{V})_k$ is not completely described by just giving that permutation. The type of $\bigwedge^m(^w\mathbf{V})_k$ (which is of course determined by $\tau$, the type of $(^w\mathbf{V})_k$) also matters. If $(^w\mathbf{V})_k$ is a full F-zip, the $k$-dimension of $\gr_C^r\bigwedge^m\mathbf{V}$ is the number of possibilities to write $r$ as a sum of $m$ distinct elements of $\supp(\tau)$ (irrespective of the order of the summands). For general F-zips, this gets even more complicated, so we leave it  at that.

\section{Dimension of \texorpdfstring{$\Hom(\1(r), (^w\mathbf{V})_k)$}{Hom(1(r), ...)}}

Notation: Fix  $r\in\Z$, a (not necessarily full) type of F-zips $\tau$ and $w\in S_n$, where $n:=\sum\tau(i)$. Set
\begin{align*}
 \mathbf{W}&:=(^w\mathbf{V})_k, \\
  T &:= C^r(\mathbf{W})\cap D_r(\mathbf{W}), \\
  L &:= \gr_C^r\mathbf{W}, \\
  R &:= \gr_r^D\mathbf{W}.
\end{align*}

Write $w=c_1\dotsm c_s$ as a product of disjoint cycles, where $c_i=(j_{i,1}\quad \dotsb \quad j_{i,m_i})$ ($m_i=\ord(c_i)$), such that $\{1,\dotsc,n\}$ is the set of all $j_{i,l}.$

\begin{Proposition}\label{prop:zahlen}
  Let $k$ be algebraically closed. There exists $d'\in \N_0$ such that
  \begin{equation*}
    \Hom(\1(r), (^w\mathbf{V})_k)\cong \bigoplus_{i\in J}\F_{p^{m_i}}\oplus k^{d'}
  \end{equation*}
  as $\F_p$-vector spaces, where
  \begin{equation*}
    J =\{i\in\{1,\dotsc,s\}\suchthat \forall \mu\in\{1,\dotsc,m_i\}\colon e_{ j_{i,\mu}}\in L\}
  \end{equation*}
  (Here $\left\{e_i\right\}_{1\leq i\leq n}$ is the basis we used to define the F-zip associated to $w$, a basis that is adapted to both $C^\bullet$ and $D_\bullet$.)

  Set
  \begin{align*}
    \psi_1(r,\mathbf{W})&:=\psi_1(r,\tau,w):=d=\sum_{i\in J} m_i, \\
    \psi_2(r,\mathbf{W})&:=\psi_1(r,\tau,w):=d'.
  \end{align*}
\end{Proposition}

\begin{Proof}
  In the following,  ``$i\in T$'' will sometimes be used as an abbreviation for ``$e_i\in T$'', etc.

  Recall that the isomorphisms $C^\nu/C^{\nu+1}\to D_\nu/D_{\nu-1}$ are given by $e_\mu\mapsto e_{w(\mu)}$, in particular
  \begin{equation*}
    i\in L\iff w(i)\in R.
  \end{equation*}

  Let $f\colon k\to k^n$ be $k$-linear.

  Write $f(1)=\sum a_{i}e_{i}.$ Here $f\in\Hom(\1(r), \mathbf{W})\implies\forall i\notin T\colon a_{i}=0.$

  We then have
  \begin{align*}
    ((e_{i}\mapsto e_{w(i)}) \circ \gr_C^rf)(1) &= \sum b_{i}e_{i}, \\
    \gr_r^Df(1) &= \sum c_{i}e_{i},
  \end{align*}
  where
  \begin{align*}
    b_{i}&=
    \begin{cases}
      a_{w^{-1}(i)}^p, & w^{-1}(i)\in L \\
      0, & w^{-1}(i)\notin L
    \end{cases} \\
    c_{i}&=
    \begin{cases}
      a_{i}, & i\in R \\
      0, & i\notin R
    \end{cases}
  \end{align*}

  So, $f\in\Hom(\1(r), \mathbf{W})$ if and only if $\left(\forall i\notin T\colon a_{i}=0 \text{ and }\sum b_{i}e_{i}=\sum c_{i}e_{i}\right)$ if and only if for all $i\in \{1,\dotsc,n\}$ the following conditions hold:
  \begin{equation*}
  \begin{cases}
    a_{i}=0, & i \notin T \\
    a_{i}=a_{w^{-1}(i)}^p, & w^{-1}(i)\in L\; (\Leftrightarrow i\in R)\\
    \text{(no extra condition)}, & i\in T, w^{-1}(i)\notin L
  \end{cases}
\end{equation*}

If $a_i$ is to be an ``$\F_{p^m}$-parameter'', we need to have a cycle of conditions
\begin{equation*}
  a_i=a_{w^{-1}(i)}^p=a_{w^{-2}(i)}^{p^2}
  =\dotsb=a_{w^{-m}(i)}^{p^m}=a_i^{p^m}
\end{equation*}

Then $w^{l}(i)\in L$ for all $l$ (hence $w^{l}(i)\in R$ for all $l$, hence $w^{l}(i)\in T$ for all $l$) and we can choose $a_i\in \F_{p^m}$ freely (provided $m$ is minimal with $w^{-m}(i)=i$) and that will determine all $a_{w^{l}(i)}$ uniquely.

If $a_i$ is not  involved in such a cycle, then either we need to enforce $a_i=0$ or $a_i$ might be in a chain of conditions, which is not circular. Such a chain will yield one ``$k$-parameter.''

Parametrizing $\Hom(\1(r),\mathbf{W})$ like that,  we may write ${\Hom(\1(r),\mathbf{W}))\cong \bigoplus_{i\in J}\F_{p^{m_i}}\oplus k^{d'}}$.

Here $d'$ is the number of $i\in\{1,\dotsc,n\}$ such that there is a non-circular chain of conditions $a_i=\dotsb=a_{w^{-m}(i)}^{p^m}$, $m\geq 0$, which can't be extended to the left or right and where the conditions above force no $a_i$ to be zero. We have such a chain if and only if:
\begin{compactitem}
\item $w^l(i)\notin T$ for some $l$,
\item $w^{-l}(i)\in T$ for all $0\leq l\leq m$,
\item $w^{-l}(i)\in L$ for all $1\leq l\leq m$
\item $w^{-(m+1)}(i)\notin L$,
\item $i\notin L$ or $w(i)\in T$.
\end{compactitem}
\end{Proof}

\section{Example}
\label{sec:example}

We consider F-zips of type $\tau$ over an algebraically closed field $k$, where $\tau(1)=3$, $\tau(0)=2$ and $\tau(i)=0$ for all $i\notin\{0,1\}$. Hence, with the notation used in the classification of F-zips, we have $(i_1,i_2)=(1,0)$, $m_1=3$, $m_2=5$.

Let $w\in S_5$. The filtrations of $(^w\mathbf{V})_k$ are such that $C^0/C^1=\langle e_4,e_5\rangle$, $C^1/C^2=\langle  e_1,e_2,e_3\rangle$, $D_0/D_{-1}=\langle e_{w(4)},e_{w(5)}\rangle$, $D_1/D_{0}=\langle e_{w(1)},e_{w(2)},e_{w(3)}\rangle$ for a basis $\left\{e_i\right\}_i$ of the underlying vector space.

In the following, we will use $i_1\dotsb i_l$ as a shorthand for $e_{i_1}\wedge\dotsb\wedge e_{i_l}$.

The filtrations of $\bigwedge^2(^w\mathbf{V})_k$ are such that $C^0/C^1=\langle 45\rangle$, $C^1/C^2=\langle 24, 15, 14, 35, 34, 25\rangle$, $C^2/C^3=\langle 13,12,23\rangle$, $D_0/D_{-1}=\langle w(4)w(5)\rangle$, etc.

The filtrations of $\bigwedge^3(^w\mathbf{V})_k$ are such that $C^1/C^2=\langle 345,145,245\rangle$, $C^2/C^3\langle 135, 235, 234, 125, 124, 134\rangle$, $C^3/C^4=\langle 123\rangle$.

The filtrations of $\bigwedge^4(^w\mathbf{V})_k$ are such that $C^2/C^3=\langle 1345, 1245,2345\rangle$ and $C^3/C^4=\langle 1234,1235\rangle$.

Now, according to the classification of F-zips, the F-zips associated to $w_1=
\begin{pmatrix}
  1 & 2 & 3 & 4 & 5 \\
  3 & 4 & 5 & 1 & 2
 \end{pmatrix}$ and $w_2=\begin{pmatrix}
  1 & 2 & 3 & 4 & 5 \\
  2 & 3 & 5 & 1 & 4
\end{pmatrix}$
respectively are not isomorphic.

However, as a calculation using Proposition \ref{prop:zahlen} shows, $\psi_1(r,\bigwedge^m(^{w_1}\mathbf{V})_k)=\psi_1(r,\bigwedge^m(^{w_2}\mathbf{V})_k)$ for  all $r,m$.

Also $\psi_2(r,\bigwedge^m(^{w_1}\mathbf{V})_k)=\psi_2(r,\bigwedge^m(^{w_2}\mathbf{V})_k)$ for  all $r,m$.

\appendix

\section{Appendix: Code}
\label{sec:appendix}

The following Python code can be used to compute the numbers $d:=\sum_{i\in J}m_i$ and $d'$ in Proposition \ref{prop:zahlen} for a given type $\tau$.

\begin{lstlisting}[language=Python,caption=fzip.py]
#!/usr/bin/env python
# -*- coding: utf-8 -*-

from itertools import permutations, product, combinations
import math

from fziputil import *

# ----------------------
tau = ((1,3),(0,2))   # format (i_j, n_j) with i_j > i_{j+1}, n_j \ne 0
# ----------------------




(m,n) = calc_mn(tau)
r = len(tau)
W = permutations(list(range(1,n+1)))


w0J = []
for i in range(1,r+1):
    for j in range(1+m[i-1], 1+m[i]):
        w0J.append(n - m[i] - m[i-1] + j)


J = []
for i in range(1,r+1):
    J.append([])
    P = permutations(list(range(m[i-1]+1, m[i]+1)))
    for p in P:
        l = []
        for j in range(1, n+1):
            if j >= m[i-1]+1 and j <= m[i]:
                l.append(p[j - (m[i-1]+1)])
            else:
                l.append(j)
        J[i-1].append(l)


for i in range(1, r):
    J[0] = [mult(x,y) for (x,y) in product(J[0], J[i])]

WJ = J[0]
del(J)

JW = cosets(WJ, W)


reps = list(map(lambda w: mult(w, w0J), JW))  # the permutations w^\circle



    

fp_dims = dict()
k_dims = dict()
for e in range(1, n+1):  # compute e'th exterior power
    for w in reps:
        sw = tuple(w)
        try:
            fp_dims[sw].append([])
            k_dims[sw].append([])
        except KeyError:
            fp_dims[sw] = [[]]
            k_dims[sw] = [[]]
    P = list(combinations(list(range(1,n+1)), e))  # we use this to map between the set of e-subsets of {1,...,n} and {1,...,n choose e}
    C = []
    D = []
    for k in range(tau[-1][0], n * tau[0][0] + 1):  # Don't calculate numbers that are guaranteed to be zero anyway.
        # Since nothing happens below tau[-1][0], this is index 0 for data structures.
        i = k - tau[-1][0]    # Hom(1(k), Lambda^e(...)) is what we will be looking at.
        C.append(set())
        D.append({})
        for I in product(tau, repeat=e):
            if sum(map(lambda x: x[0], I)) != k:
                continue
            # we want (i_{j_1},...,i_{j_e}) with i_{j_1}+...+i_{j_e}=k
            C[i] |= build_C_elements(I, tau)
            for w in reps:
                sw = tuple(w)
                try:
                    D[i][sw] |= build_D_elements(I, tau, w)
                except KeyError:
                    D[i][sw] = build_D_elements(I, tau,  w)
       
    for k in range(tau[-1][0], n * tau[0][0] + 1):
        i = k - tau[-1][0]
        for w in reps:
            sw = tuple(w)
            fp_dims[sw][e-1].append(0)
            k_dims[sw][e-1].append(0)
        if C[i] and D[i]:
            Chigher = set()
            try:
                Chigher = C[i+1]
            except IndexError:
                pass
            except KeyError:
                pass
            L = C[i] - Chigher
            Ln = set(P.index(tuple(x))+1 for x in L)  # mapping it back to plain numbers using P
            for w in reps:
                sw = tuple(w)
                k_dim, fp_dim = 0, 0
                Dlower = set()
                try:
                    Dlower = D[i-1][sw]
                except IndexError:
                    pass
                except KeyError:
                    pass
                T = C[i] & D[i][sw]
                Tn = set(P.index(tuple(x))+1 for x in T)
                R = D[i][sw] - Dlower
                Rn = set(P.index(tuple(x))+1 for x in R)
                assert(len(Rn)==len(Ln))
                Z = Ln
                q = iota(P,e, w)
                cyc = cycles(q)
                lcyc = list(cyc)
                fp = set(range(1,nCr(n,e)+1)) - set(flatten(lcyc[:]))
                for x in fp:
                    if x in Z:
                        fp_dim += 1
                for c in lcyc:
                    success = True
                    for x in c:
                        if x not in Z:
                            success = False
                            break
                    if success:
                        fp_dim += len(c)
                    else:
                        lc = len(c)
                        for j in range(lc):
                            if c[j] not in Tn:
                                continue
                            if c[j] in Ln:
                                continue
                            m = 0
                            while c[(j-m-1) % lc] in Tn and c[(j-m-1) % lc] in Ln:
                                m += 1
                            if c[(j-m-1) % lc] in Ln:
                                continue
                            k_dim += 1

                fp_dims[sw][e-1][i] = fp_dim
                k_dims[sw][e-1][i] = k_dim


print("\n\nk")
for sw in k_dims:
    print(sw, k_dims[sw])
print("\n\nF_p")
for sw in fp_dims:
    print(sw, fp_dims[sw])
print("\n")




s = set()
for v in fp_dims.values():
    l = [tuple(i) for i in v]
    s.add(tuple(l))

\end{lstlisting}

\begin{lstlisting}[language=Python,caption=fziputil.py]
import sys, math



def mult(l1, l2):
    """multiply two permutations"""
    return [l1[l2[i]-1] for i in range(len(l1))]


def coset(subgrp, p):
    """length-minimal representative of the right coset subgrp*p"""
    s = set()
    for q in subgrp:
        s.add(tuple(mult(q,p)))
    minrep = None
    minlen = sys.maxsize
    for p in s:
        l = length(p)
        if l < minlen:
            minlen = l
            minrep = p
    return minrep


def cosets(subgrp, grp):
    """length-minimal representatives of the right cosets in subgrp\grp"""
    return {tuple(coset(subgrp, p)) for p in grp}


def length(p):
    """length of a permutation"""
    l = 0
    n = len(p)
    for i in range(n):
        for j in range(i+1,n):
            if p[i] > p[j]:
                l +=1
    return l


#  from <http://lists.canonical.org/pipermail/kragen-hacks/2013-August/000560.html>
def cycles(p):
    """Analyze permutation into a product of disjoint cycles.

    This is very straightforward to do in linear time using Python
    set objects, but using Python set objects to find the start of
    each new cycle produces the cycles in unpredictable order.
    The algorithm here should still be linear-time (ii visits each
    index at most twice; min_leftover visits each index once; the
    initial leftovers set is produced in linear time; every other
    simple statement in the function is O(1)) and also produce the
    nonempty cycles in a deterministic order, with a deterministic
    starting point.

    """
    leftovers = set(p)
    min_leftover = 1

    while leftovers:
        while min_leftover not in leftovers:
            min_leftover += 1
            assert min_leftover <= len(p)

        ii = min_leftover
        cycle = []
        cycle_set = set()
        while ii not in cycle_set:
            cycle.append(ii)
            cycle_set.add(ii)
            leftovers.remove(ii)
            ii = p[ii-1]

        if len(cycle) > 1:
            yield cycle


def calc_mn(tau):
    """given a type tau, determine its height and the numbers (m_i)"""
    m = [0]
    for (mu, t) in zip(m, tau):
        m.append(mu + t[1])
    n = sum(map(lambda x: x[1], tau))
    return (m,n)


def build_C_elements(I, tau):
    """set of all {j_1,...,j_e} where e_{j_l} in C^{i_l}, I=(i_1,...,i_e)"""
    if not I:
        return set([frozenset()])
    i = tau.index(I[0])
    tail = build_C_elements(I[1:], tau)
    L = set()
    (m,n) = calc_mn(tau)
    for j in range(m[i+1]):
        for s in tail:
            if j+1 not in s:
                t = set(s)
                t.add(j+1)
                L.add(frozenset(t))
    return L


def build_D_elements(I, tau, w):
    """set of all {j_1,...,j_e} where e_{j_l} in D_{i_l} (which depends on tau), I=(i_1,...,i_e)"""
    if not I:
        return set([frozenset()])
    i = tau.index(I[0])
    tail = build_D_elements(I[1:], tau, w)
    L=set()
    m = [0]
    (m,n) = calc_mn(tau)
    for j in range(n-m[i]):
        for s in tail:
            if w[m[i]+j] not in s:
                t = set(s)
                t.add(w[m[i]+j])
                L.add(frozenset(t))
    return L


def flatten(lst):
    """flatten a list"""
    for elem in lst:
        if type(elem) in (tuple, list):
            for i in flatten(elem):
                yield i
        else:
            yield elem


def iota(P,e,p):
    """the natural action of p in S_n on e-subsets of {1,...,n},
    P is used to map between the set of e-subsets of {1,...,n} and {1,...,n choose e}"""
    n = len(p)
    q = []
    if e == 0:
        return q
    for S in P:
        T = tuple(sorted([p[j-1] for j in S]))
        q.append(P.index(T)+1)
    return q


def nCr(n,r):
    """binomial coefficient"""
    f = math.factorial
    return int(f(n) / f(r) / f(n-r) + 0.01)
\end{lstlisting}

\printbibliography

\end{document}